\documentclass[11pt,letterpaper]{article}
\usepackage[utf8]{inputenc}
\usepackage[T1]{fontenc}
\usepackage[a4paper, left=5cm, textwidth=11.0cm,  textheight=22.2cm,
 showframe=false,verbose=true]{geometry}
\usepackage{marginnote}
\RequirePackage{tocloft}

\RequirePackage{tabularx,booktabs}
\usepackage[fleqn]{amsmath}
\usepackage{amssymb}
\usepackage{amsthm}
\usepackage{mathtools}
\usepackage{listings}
\usepackage{float}
\usepackage[hyphens]{url}
\usepackage{caption}
\usepackage{subfig}
\usepackage[verbose]{placeins}
\usepackage{fancyvrb}
\usepackage[labelformat=empty]{caption}
\allowdisplaybreaks
\renewcommand{\Im}{\mathop{\rm Im}\nolimits}
\renewcommand{\Re}{\mathop{\rm Re}\nolimits}

\newcommand{\Ee}{ \mathrm{e} }
\DeclareMathOperator{\sech}{sech}
\RequirePackage[dvipsnames,table]{xcolor}
\definecolor{formula}{rgb}{.6,0,0}
\def\OEIS{\textsc{\small{{OEIS}}}~}
\usepackage{relsize}
\def\Kappa{\mathlarger{\mathlarger{\mathlarger{\kappa}}}}
\newcommand{\seqnum}[1]{\href{http://www.oeis.org/#1}{\small#1}}

\usepackage{imakeidx}
\usepackage{xparse}
\NewDocumentCommand{\stirling}{omm}
{%
    \genfrac\{\}{0pt}{}{#2}{#3}%
    \IfValueT{#1}{_{\!#1}}%
}%
\usepackage{graphicx}
\usepackage{tikz}

\usepackage{epstopdf}
\RequirePackage[pdftex,hyperfootnotes=false,pdfpagelabels]{hyperref}
\pdfcompresslevel=9 \pdfadjustspacing=1 \hypersetup{colorlinks=true,linktoc=all,%
    urlcolor=Blue,linkcolor=Blue,citecolor=Blue,linktocpage=true,pdfstartpage=1,pdfstartview=FitV,%
    breaklinks=true,pdfpagemode=UseNone,pageanchor=true,pdfpagemode=UseOutlines,plainpages=false,%
    bookmarksnumbered,bookmarksopen=true,bookmarksopenlevel=0,hypertexnames=true,pdfhighlight=/O,%
    pdftitle={An introduction to the Bernoulli function},%
    pdfauthor={Peter H. N. Luschny},pdfsubject={Exploring a variant of the zeta function interpolating the Bernoulli numbers based on an integral representation suggested by J. Jensen.},%
    pdfkeywords={Bernoulli function, Bernoulli constants, Bernoulli numbers, Bernoulli cumulants, Bernoulli functional equation, Bernoulli central function, extended Bernoulli function, Bernoulli central polynomials, alternating Bernoulli function, Stieltjes constants, Riemann zeta function, Hurwitz zeta function, Worpitzky numbers, Worpitzky transform, Hasse representation, Hadamard product,  Genocchi numbers, Genocchi function, Euler secant numbers, Euler tangent numbers, Euler zeta numbers, Euler function, Euler tangent function, Euler secant function, Eulerian numbers, Andr\'e numbers, Andr\'e function, Seki numbers, Seki function, Swiss-knife polynomials},%
    pdfcreator={pdfLaTeX}}
\usepackage[all]{hypcap}

\normalsize

\makeatletter
\newcommand{\mathleft}{\@fleqntrue\@mathmargin0pt}
\newcommand{\mathcenter}{\@fleqnfalse}
\makeatother
\usepackage{accents}

\usepackage{etoolbox}
\usepackage{comment}
\usepackage{setspace}
\usepackage{microtype}
\usepackage{rotating}
\usepackage{makeidx}\makeindex[title=Index]

\usepackage{tocloft}
\makeatletter
\renewcommand{\@cftmaketoctitle}{}
\makeatother
\makeindex[columns=2, title=Index]
\newcounter{mysection}
\newcommand{\sect}[2]{
    \stepcounter{mysection}
    \vspace*{#2 pt}  {\noindent\large\textbf{\textcolor{MidnightBlue}{\arabic{mysection} -- #1}}}
    \vspace*{12pt}    \phantomsection\addcontentsline{toc}{section}{\arabic{mysection} -- #1}}
\usepackage{datetime}
\usepackage{layout}

\begin{document} 
\title{{\textls[16]
            { \textsc{\textbf{\textcolor{MidnightBlue}{An Introduction to the \\ Bernoulli Function}}}}}}
\author{\textsc{Peter H. N. Luschny}}
\date{}
\maketitle
\vspace{-18pt}
\begin{spacing}{0.90}
    \begin{center} \parbox{266pt}{ {\footnotesize\noindent \textcolor{MidnightBlue}{\textsc{Abstract. }}
         We explore a variant of the zeta function interpolating the Bernoulli numbers based on an integral 
         representation suggested by J. Jensen. The \textit{Bernoulli function} \(\operatorname{B}(s, v) = - s\, \zeta(1-s, v)\) 
         can be introduced independently of the zeta function if it is based on a formula first
         given by Jensen in 1895.  We examine the functional equation of \(\operatorname{B}(s, v) \) and
        its representation by the Riemann \(\zeta\) and \(\xi\) function, and
        recast classical results of Hadamard, Worpitzky, and Hasse in terms of \(\operatorname{B}(s, v).\)
        The \textit{extended Bernoulli function} defines the Bernoulli numbers
        for odd indices basing them on rational numbers studied by Euler in 1735
        that underlie the \textit{Euler and Andr\'{e} numbers}.
        The \textit{Euler function} is introduced as the difference between values of the
        Hurwitz-Bernoulli function.  The \textit{Andr\'{e} function} and the \textit{Seki function}
        are the unsigned versions of the extended Euler resp. Bernoulli function.
        }}
    \end{center}
\end{spacing}

\vspace*{16pt}\hspace*{18pt}
\textcolor{Bittersweet}{
    \textsc{Notations}~\pageref{page-synopsis}
    \textsc{Index}~\pageref{page-index}
    \textsc{Plots}~\pageref{page-plots}
    \textsc{Contents}~\pageref{page-contents}
}

\bigskip

\sect{Prologue: extension by interpolation}{20}

\noindent \textit{\textbf{\textcolor{MidnightBlue}{The question}}}
\bigskip

Andr\'e Weil recounts the origin of the gamma function
in his historical exposition of  number theory \cite[p.\ 275]{Weil}:
\begin{quote}
    \textit{``Ever since his early days in Petersburg Euler had been interested in the interpolation
        of functions and formulas given at first only for integral values of the argument;
        that is how he had created the theory of the gamma function.''}
\end{quote}
The three hundred year success story of Euler's gamma function shows how
fruitful this question is. The usefulness of such an investigation is not limited
by the fact that there are infinitely many ways to interpolate a sequence of numbers.

\smallskip
This essay will explore the question: how can the Bernoulli numbers be interpolated most meaningfully?

\newpage
\noindent \textit{\textbf{\textcolor{MidnightBlue}{The method}}}

\bigskip

The Bernoulli numbers had been known for some time at the beginning of the 18th century
and used in the (now Euler--Maclaurin called) summation formula in analysis, first
without realizing that these numbers are the same each case.

\smallskip
Then, in 1755, Euler baptized these numbers \textit{Bernoulli numbers}
in his \textit{Institutiones calculi differentialis}
(following the lead of de Moivre).
After that, things changed, as Sandifer \cite{Sandifer} explains:

\begin{quote}
    \textit{``... for once the Bernoulli numbers had a name,
        their diverse occurrences could be recognized, organized, manipulated and
        understood. Having a name, they made sense.''}
\end{quote}

The function we are going to talk about is not new.
However, it is not treated as a function in its own right and with its own name.
So we will give the beast a name. We will call the interpolating function the
\textit{Bernoulli function}. Since, as Mazur \cite{Mazur} asserts,
the Bernoulli numbers ``act as a unifying force, holding together seemingly
disparate fields of mathematics,'' this should all the more
be manifest in this function.

\bigskip

\noindent \textit{\textbf{\textcolor{MidnightBlue}{What to expect}}}

\bigskip

This note is best read as an annotated formula collection; we refer to the cited literature for the proofs.

\smallskip
The Bernoulli function and the Riemann zeta function can be understood as complementary pair.
One can derive the properties of one from those of the other. For instance, all the questions Riemann
associated with the zeta function can just as easily be discussed with the Bernoulli function.
In addition, the introduction of the Bernoulli function leads to a better understanding which
numbers the \textit{generatingfunctionologists} \cite{Wilf} should hang up on their clothesline
(see the discussion in \cite{Luschny}).

\smallskip
Seemingly the first to treat the Bernoulli function in our sense was Jensen \cite{Jensen},
who gave an integral formula for the Bernoulli function of remarkable simplicity.
Except for referencing Cauchy's theorem, he did not develop the proof. The proof
is worked out in Johansson and Blagouchine \cite{Johansson}.

\smallskip

A table of {contents}~(\pageref{page-contents}) can be found at the end 
of the paper.
\clearpage \newpage 
\hspace{-24pt}{ \sect{Notations and jump-table to the definitions}{0}
\label{page-synopsis}
 \begingroup
 \hypersetup{linkcolor=blue}
\vspace{-18pt}
{\setlength{\mathindent}{0pt}
    \\ \marginnote{\textit{eq.no.} }
    \begin{equation*}
        \operatorname{B}(s,v) \, = \,  2 \pi \int_{- \infty}^\infty  \frac{(v-1/2 + i z)^{s}} {({\Ee}^{- \pi z} + {\Ee}^{\pi z})^2}\, dz
        \hspace{67pt} \text{Jensen formula}
        \nonumber \marginnote{\(\rightarrow\) \underline{\ref{int5}}}
    \end{equation*}
    \begin{equation*}
        \operatorname{B}(s) \, = \, \operatorname{B}(s,1) \, = \, -s\, \zeta(1-s)
        \hspace{84pt} \text{Bernoulli function}
        \marginnote{\(\rightarrow\)\phantom{8} \underline{\ref{Bdef}}}
    \end{equation*}
    \begin{equation*}
        \operatorname{B}^c(s) \, = \, \operatorname{B}\left( s,1/2 \right)
        \hspace{97pt} \text{Central Bernoulli function}
        \marginnote{\(\rightarrow\) \underline{\ref{cbrep}}}
    \end{equation*}
    \begin{equation*}
        \operatorname{B}^{*}(s) \, = \, \operatorname{B}(s)(1-2^s)
        \hspace{70pt} \text{Alternating Bernoulli function}
        \marginnote{\(\rightarrow\) \underline{\ref{betastar}}}
    \end{equation*}
    \begin{equation*}
        \operatorname{B}_n \, = \, \operatorname{B}(n)
        \hspace{171pt} \text{Bernoulli numbers}
        \marginnote{\(\rightarrow\)\phantom{8} \underline{\ref{bernum}}}
    \end{equation*}
    \begin{equation*}
        \operatorname{B}_{\sigma}(s)  =  \frac{2^{s - 1}}{2^s - 1} \left(\operatorname{B}\left(s, \tfrac34 \right) - \operatorname{B}\left(s, \tfrac14 \right) \right)  
        \hspace{22pt} \text{Bernoulli secant function}
        \marginnote{\(\rightarrow\) \underline{\ref{bernsec}}}
    \end{equation*}
    \begin{equation*}
        \operatorname{G}(s,v) \, = \,  2^s \left( \operatorname{B}\left(s,\tfrac{v}{2}\right) - \operatorname{B}\left(s,\tfrac{v+1}{2}\right) \right)
        \hspace{28pt} \text{Gen. Genocchi function}
        \marginnote{\(\rightarrow\) \underline{\ref{GG}}}
    \end{equation*}
    \begin{equation*}
        \operatorname{G}(s) \, = \, \operatorname{G}(s,1)\,=\,  2 (1 - 2^s) \operatorname{B}(s)
        \hspace{69pt} \text{Genocchi function}
        \marginnote{\(\rightarrow\) \underline{\ref{GenocciFun2}}}
    \end{equation*}
    \begin{equation*}
        \operatorname{E}(s,v) \, = \, - \frac{\operatorname{G}(s+1,v)}{s+1}
        \hspace{73pt}  \text{Generalized Euler function}
        \marginnote{\(\rightarrow\) \underline{\ref{EG}}}
    \end{equation*}
    \begin{equation*}
        \operatorname{E}_{\tau}(s) \, = \,  \operatorname{E}(s,1) 
        \hspace{129pt}  \text{Euler tangent function}
        \marginnote{\(\rightarrow\) \underline{\ref{eulsimpl}}}
    \end{equation*}
    \begin{equation*}
        \operatorname{E}_{\sigma}(s) \, = \,  2^{s} \operatorname{E}(s, \tfrac{1}{2}) 
        \hspace{122pt}  \text{Euler secant function}
        \marginnote{\(\rightarrow\) \underline{\ref{EulerCentral}}}
    \end{equation*}
    \begin{equation*}
        \operatorname{E}_n \, = \,  \operatorname{E}_{\sigma}(n)  
        \hspace{184pt}  \text{Euler numbers}
        \marginnote{\(\rightarrow\) \underline{\ref{EulerNum}}}
    \end{equation*}
    \begin{equation*}
        \widetilde{\zeta}(s)  \, = \,   \zeta \left({s} \right) + \frac{ \zeta \left({s},{\frac{1}{4}} \right)  - \zeta \left({s},{\frac{3}{4}} \right) } {2 ^{s} - 2}
        \hspace{52pt} \text{Extended zeta function}
        \nonumber \marginnote{\(\rightarrow\) \underline{\ref{lzeta}}}
    \end{equation*}
    \begin{equation*}
        \mathcal{B}(s) \, = \, -s \, \widetilde{\zeta}(1-s)
        \hspace{91pt} \text{Extended Bernoulli function}
        \nonumber \marginnote{\(\rightarrow\) \underline{\ref{genbern}}}
    \end{equation*}
    \begin{equation*}
        \mathcal{E}(s) \, = \,  \frac{4^{s+1} - 2^{s+1}}{ (s+1)!} \mathcal{B}(s+1)
        \hspace{66pt}  \text{Extended Euler function}
        \nonumber \marginnote{\(\rightarrow\) \underline{\ref{exteulbern}}}
    \end{equation*}
    \begin{equation*}
        \mathcal{A} (s) \, = \, (-i)^{s+1} \operatorname{Li}_{-s}(i) \, +\,  i^{s+1} \operatorname{Li}_{-s}(-i)  
        \hspace{56pt}  \text{André function}
        \nonumber \marginnote{\(\rightarrow\) \underline{\ref{UnsignedAndreFun}}}
    \end{equation*}
\begin{equation*} 
    \mathcal{A}^{\ast}(s) \, =\, i e^{-\frac{1}{2} i \pi  s} \left(\text{Li}_{-s}(-i)-e^{i \pi  s} \text{Li}_{-s}(i)\right)  
     \hspace{7pt}  \text{Signed André function}
    \nonumber \marginnote{\(\rightarrow\) \underline{\ref{SignedAndreFun}}}
\end{equation*} 
   \begin{equation*}
        \mathcal{S} (s) \, =  \frac{s \, \mathcal{A} (s-1)}{4^s - 2^s}  
        \hspace{162pt} \text{Seki function}
        \marginnote{\(\rightarrow\) \underline{\ref{SekiAndreFun}}} 
    \end{equation*}
\begin{equation*} 
    \mathcal{S}^{\ast}(s) \, =\, \frac{s\, e^{i \pi s / 2}}{2^s - 4^s}  \left( e^{i \pi s} \text{Li}_{1 - s}(-i) + \text{Li}_{1 - s}(i) \right)
    \hspace{8pt} \text{Signed Seki function}
    \marginnote{\(\rightarrow\) \underline{\ref{SignedSekiFun}}}
\end{equation*} 
    \begin{equation*}
    \mathcal{Z}_n \, = \, \frac{\mathcal{A}_n}{n !} 
    \hspace{172pt}  \text{Euler zeta numbers}
    \nonumber \marginnote{\(\rightarrow\) \underline{\ref{eulzetanum}}}
\end{equation*}
    \begin{equation*}
        \operatorname{\xi}(s) \, = \, \operatorname{B}(s)\, \frac{\sigma !}{\pi^{\sigma}}  \ \text{with } \sigma =  (1-s)/2
        \hspace{54pt} \text{Riemann }\xi\text{ function}
        \nonumber \marginnote{\(\rightarrow\) \underline{\ref{q554}}}
    \end{equation*}
    \begin{equation*}
        \gamma \, = \,  - \operatorname{B}^{'}(0)  \  \text{(Euler),} \quad
        \gamma_{n-1} \, = \,  - \frac{1}{n}\operatorname{B}^{(n)}(0) \quad
        \hspace{11pt} \text{Stieltjes constants} \nonumber
        \marginnote{\(\rightarrow\) \underline{\ref{q26}}}
    \end{equation*}
    \enlargethispage*{20pt}
   }  \endgroup
 \clearpage \newpage

\!\begin{figure}
    \centering \includegraphics[width=150pt]{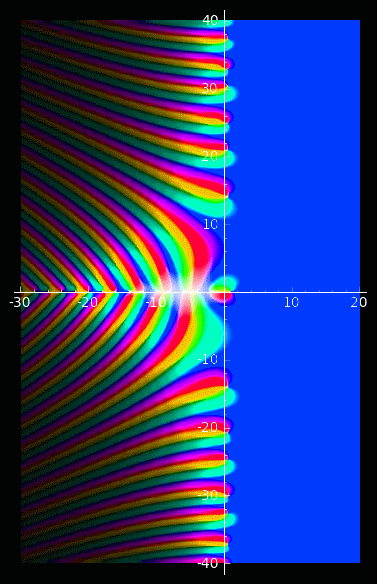}%
    \( \)
    \centering \includegraphics[width=150pt]{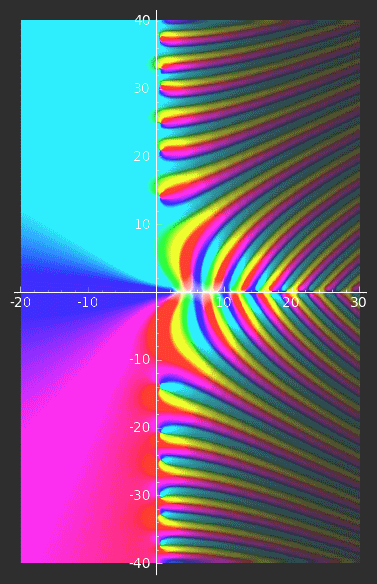}%
    \caption[Zeta versus Bernoulli function]{\textit{The zeta and the Bernoulli function in the complex plane.}}%
\end{figure}

\sect{Stieltjes constants and Zeta function}{0}

The generalized Euler constants, also called \textit{Stieltjes constants}\index{Stieltjes!constants},
are the real numbers \( \gamma_n \) defined by the Laurent series in a neighborhood of \(s=1\)
of the Riemann zeta function\index{Riemann!\(\zeta\)-function} (\cite{DLMF, RiemannEng}),
\begin{equation}
    \zeta(s)= \frac{1}{s-1}+\sum_{n=0}^\infty \frac{(-1)^n}{n!} \gamma_n (s-1)^n, \qquad (s \ne 1) .
    \label{q1}
\end{equation}
As a particular case they include the Euler constant\index{Constant!\(\gamma\ \)Euler}
\( \gamma_0 = \gamma \ \, ( \approx \frac{228}{395}) \).
Blagouchine \cite{BlagouchineExpansions}
gives a detailed discussion with many historical notes.
Following Franel \cite{Franel}, Blagouchine \cite{BlagouchineEval}  shows that
\begin{equation}
    \gamma_n \,=\,  - \frac{2\, \pi}{n+1} \int\limits_{- \infty}^{+\infty}\,
    \frac{ \log\!\left(\frac{1}{2} + iz \right)^{n+1}}
    {\left( {{\Ee}^{- \pi \,z}}+{{\Ee}^{\pi \,z }} \right) ^{2}}\,{dz} .
    \label{q2}
\end{equation}
Here, and in all later similar formulas, we write \( f(x)^n\) for \( (f(x))^n\)
and take the principal branch of the logarithm implicit in the exponential.

Recently Johansson  \cite{JohanssonHurZeta} observed that one can employ the integral representation
(\ref{q2}) in a particularly efficient way to approximate the Stieltjes constants with prescribed precision numerically.

\newpage

\begin{figure}
    \centering
    \includegraphics[width=310pt]{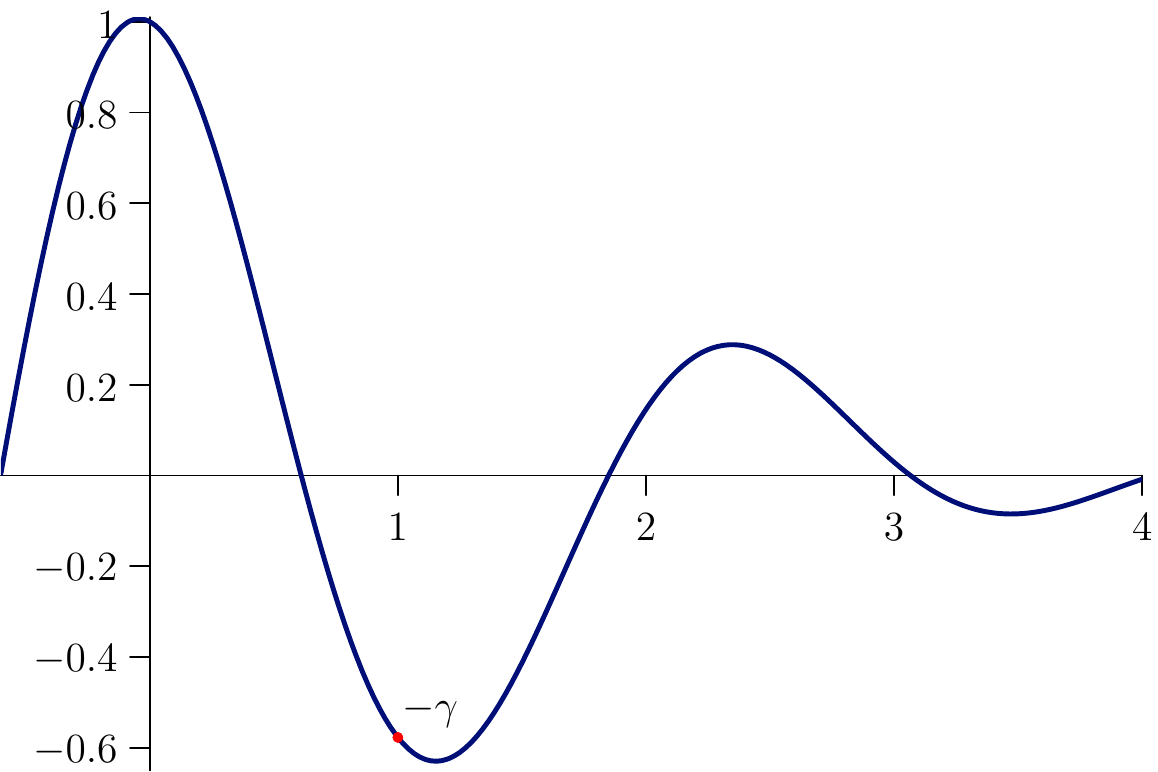}%
    \caption[Bernoulli constants as a function]{\textit{The Bernoulli constants seen as values of a real function.}}%
    \label{fig:fig2benrconst}
\end{figure}

\sect{Bernoulli constants and Bernoulli function}{0}

\noindent The \textit{Bernoulli constants}\index{Bernoulli!constants}
\( \beta_n \)  are defined for \( n \ge 0 \) as
\begin{equation} {
    \beta_n \,  = \,  2\pi \! \int\limits_{- \infty }^{+\infty }{\frac {
        \log\! \left( \frac{1}{2} + iz \right) ^n }
    {\left( {{\Ee}^{- \pi \,z}}+{{\Ee}^{\pi \,z }} \right) ^{2}}}\,{dz} . }
    \label{q3}
\end{equation}

\noindent The \textit{Bernoulli function}\index{Bernoulli!function} is defined as
\begin{equation}
    \operatorname{B}(s) \, = \, \sum _{n=0}^{\infty }\beta_{n} \frac {s^n}{n!}\,.
    \label{Bdef}
\end{equation}
Since \( \beta_0 = 1 \) we have in particular \( \operatorname{B}(0) = 1 \).

Comparing the definitions of the constants \(\gamma_n\) and \(\beta_n\),
we see that in (\ref{q3}) the exponent
of the power is synchronous to the index and the factor \( -1/(n+1)\) in (\ref{q2})
has disappeared. In other words, \( \gamma_n = - \frac{\beta_{n+1}}{n+1} \),
and  the \textit{Riemann zeta representation} of the Bernoulli function
becomes

\begin{align}
    \operatorname{B}(s)
     & =  \sum _{n=0}^{\infty }\beta_{n} \frac {s^n}{n!}  \nonumber              \\
     & = 1+  \sum _{n=1}^{\infty } \beta_{n} \frac{s^n}{n!} \nonumber            \\
     & = 1+  \sum _{n=0}^{\infty } \beta_{n+1} \frac{s^{n+1}}{(n+1)!}  \nonumber \\
     & = 1 - \sum_{n=0}^\infty \gamma_n \frac{s^{n+1}}{n!} \nonumber             \\
     & = 1 - s \sum_{n=0}^\infty \gamma_n \frac{s^n}{n!} \nonumber               \\
     & =-s\zeta(1-s) . \label{Bzeta}
\end{align}
The singularity of  \(  -s\zeta(1-s) \) at \( s=0 \) is removed by
\( \operatorname{B}(0) = 1 \). Thus \( \operatorname{B}(s) \) is an entire
function with its defining series converging everywhere in \( \mathbb{C}\).

\sect{The Bernoulli numbers}{20}

We define the \textit{Bernoulli numbers}\index{Bernoulli!numbers} as the values of the
Bernoulli function at the nonnegative integers. According to (\ref{Bdef}) this means
\begin{equation}
    \operatorname{B}_n \, = \,  \operatorname{B}(n)  \, = \, \sum _{j=0}^{\infty }\beta_{j} \frac {n^j}{j!}.
    \label{bernum}
\end{equation}
Since for \( n > 1 \) an odd integer \( -n\zeta(1-n) = 0 \), the Bernoulli numbers
vanish at these integers and
\( \operatorname{B}_{n} = 1 - n \sum_{j=0}^\infty \gamma_j \frac{n^j}{j!} \) implies%
\begin{equation}
    \sum_{j=0}^\infty \gamma_j \frac{n^j}{j!} = \frac{1}{n}   \text{ and }
    \sum_{j=0}^\infty \beta_j \frac{n^j}{j!} = 0 \qquad (n > 1 \text{ odd}).
    \label{a1}
\end{equation}

\sect{The expansion of the Bernoulli function}{20}

\noindent The Bernoulli function \( \operatorname{B}(s) = -s\zeta(1-s) \)
can be expanded by using the generalized Euler--Stieltjes constants
\begin{equation}
    \operatorname{B}(s)\, =\, 1 - \gamma s - \gamma_1 s^2 - \frac{\gamma_2}{2}
    s^3 - \frac{\gamma_3}{6} s^4 \,\ldots \  . \label{q7}
\end{equation}
Or, in its natural form, using the Bernoulli constants
\begin{equation}{
        \operatorname{B}(s)  \, = \, 1 + \beta_1 s + \frac{\beta_2}{2} s^2
        + \frac{\beta_3}{6} s^3 + \frac{\beta_4}{24} s^4 \,\ldots \  . }
    \label{q8}
\end{equation}

{\small
\begin{table}[t]
    \centering
    \begin{tabular}[c]{c|r}
        \( r \)   & \( \beta(r) \quad \qquad \qquad \qquad \qquad \qquad \)        \\ \hline
        \(- 1\)   & \texttt{-1.0967919991262275651322398023421657187190}\(\dots \) \\
        \(- 1/2\) & \texttt{0.3000952439768656513643742483305378454480}\(\dots\)   \\
        \( 0 \)   & \texttt{1.0}  \hspace{230pt}                                   \\
        \( 1/2 \) & \texttt{ 0.2364079388130323148951169845913737350793}\(\dots \) \\
        \( 1 \)   & \texttt{-0.5772156649015328606065120900824024310421}\(\dots \) \\
        \( 3/2 \) & \texttt{-0.4131520868458801199329318166967102536980}\(\dots \) \\
        \( 2 \)   & \texttt{ 0.1456316909673534497211727517498026382754}\(\dots \) \\
        \( 5/2 \) & \texttt{ 0.2654200629584708272795102903586382709016}\(\dots \) \\
        \( 3 \)   & \texttt{ 0.0290710895786169554535911581056375880771}\(\dots \) \\
        \( 7/2 \) & \texttt{-0.0845272473711484887663180676975841853310}\(\dots \) \\
        \( 4 \)   & \texttt{-0.0082153376812133834646401861710135371428}\(\dots \)
    \end{tabular}
    \caption{\textit{The Bernoulli constants for some rational \( r \).}}
\end{table}
}

Although we will often refer to the well-known properties of the zeta function
when using (\ref{Bzeta}), our definition of the Bernoulli function and the
Bernoulli numbers only depends on (\ref{q3}) and (\ref{Bdef}).

The index \( n \) in definition (\ref{q3}) is not restricted to integer values.
For illustration the function \( \beta_r \) is plotted in figure \underline{\ref{fig:fig2benrconst}},
where the index \( r \) of \( \beta \) is understood to be a real number.
The table above displays some numerical values of Bernoulli constants.

\sect{Integral formulas for the Bernoulli constants}{20}

\noindent Let us come back to the definition of \( \beta_n \) as given in (\ref{q3}).
The appearance of the imaginary unit forces complex integration; on the
other hand, only the real part of the result is used. Fortunately, the definition
can be simplified such that the computation stays in the realm of reals
provided \( n \) is a nonnegative integer.

Using the symmetry of the integrand with respect to the \( y \)-axis and
\( \left( {{\Ee}^{- \pi \,z}} + {{\Ee}^{\pi \,z }} \right) ^{2} = 4 \cosh(\pi z)^2 \)
we get from the definition (\ref{q3})
\begin{equation}
    \beta_n \ = \ \pi \! \int_{0}^{\infty}\,
    {\frac{\Re\! \left(\log(\frac12 + iz)^n \right) }{ \cosh(\pi z) ^2 }} dz \, .
    \label{q9}
\end{equation}
For the numerator of the integrand we set for \( n \ge 0 \)
\begin{equation}
    \sigma_n(z) \, = \, \Re (\log ( \tfrac{1}{2} + iz )^n ).
    \label{q10}
\end{equation}
By induction we see that
\begin{equation}
    \sigma_n(z) = \sum_{ k = 0 }^{ \lfloor n/2 \rfloor }
    (-1)^{k} \binom{n}{2k} a(z)^{n - 2k}\, b(z)^{2 k},
    \label{q11}
\end{equation}
where \( a(z) = \log(z^2+ \frac{1}{4})/2 \) and \( b(z) = \arctan(2z) \).

\begin{figure}
    \centering
    \includegraphics[width=300pt]{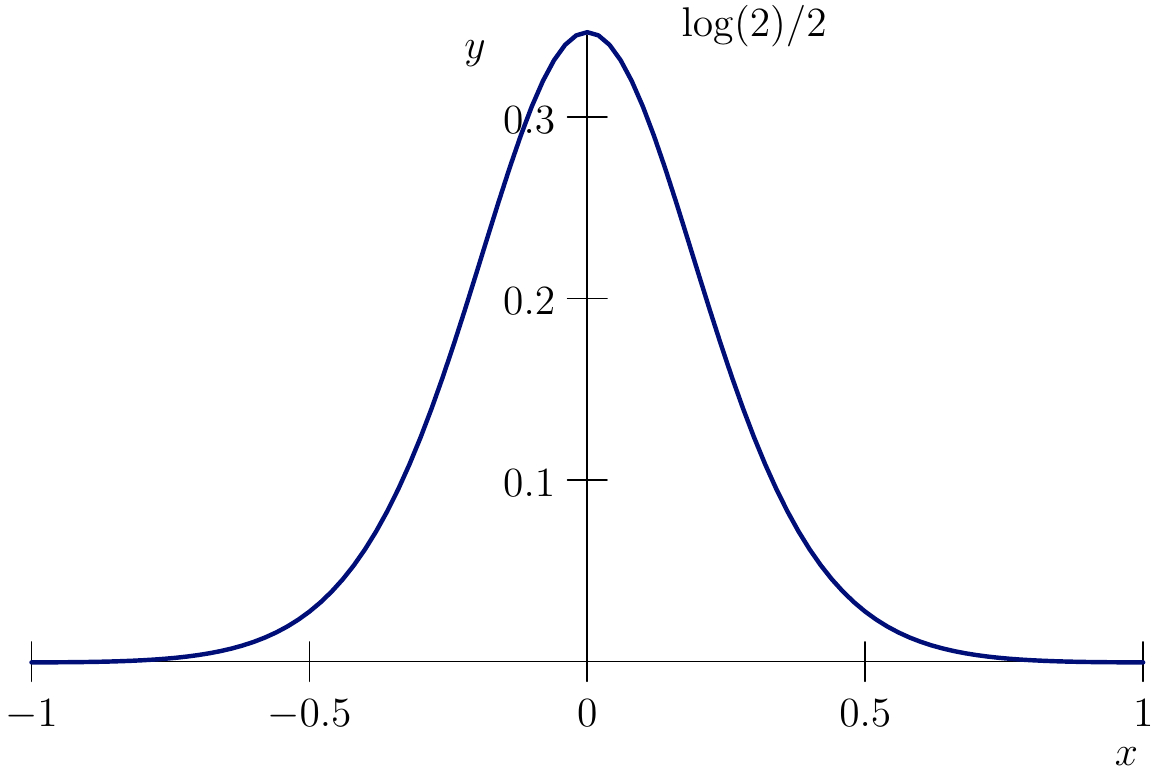}%
    \caption[Integral representation of \(\gamma/\pi\)]
    {\(f(z) = -\log((z^2+ \frac14)^{\frac14})\sech(\pi z)^2 \) \textit{with}
        \( \int_{-\infty }^{\infty} f(z)\, dz = \frac{\gamma}{ \pi}.\)}%
     \label{fig:fig3gammapi}
\end{figure}

\noindent Therefore the constants can be computed by real integrals,
\begin{equation}
    \beta_n \,=\, \pi  \! \int_{0}^{\infty} \frac{\sigma_n(z)}{\cosh(\pi z)^2}\, dz \, .
    \label{q12}
\end{equation}
Similarly we have for the Stieltjes constants
\begin{equation}
    \gamma_n \,=\, - \frac{\pi}{n+1} \int_{0}^{\infty}
    \frac{\sigma_{n+1}(z)}{\cosh(\pi z)^2}\, dz \,.
    \label{q13}
\end{equation}

\sect{Some special integral formulas}{20}

Let us consider some special cases of these integral formulas.
Formula (\ref{q13}) reads for \(n=0 \)
\begin{equation}
    \gamma \, =\, - \pi \int_{0}^{\infty} \frac{\log(z^2+ \tfrac{1}{4} ) }{2 \cosh(\pi z)^2} \, dz  .
    \label{q14}
\end{equation}
This follows since for real \(z\)
\begin{equation}
    \sigma_1(z) = \Re \log\left(\frac12+i z\right)  =
    \log\left(\frac12\sqrt{4 z^2 + 1}\right) = \frac12\log\left(z^2+\frac14\right). \nonumber
\end{equation}
Using the symmetry of the integrand with respect to the \( y \)-axis this can be
rephrased as: Euler's gamma is \( \pi \) times the integral of
\(  - \log((z^2+ \tfrac{1}{4} )^{\tfrac{1}{4}})\sech(\pi z)^2 \) over the real line.
See figure \underline{\ref{fig:fig3gammapi}}.

With the abbreviations
\( a = \log(z^2+ \frac{1}{4})/2, \, b = \arctan(2z) \) and \( c = \cosh(\pi z) \)
the first few Bernoulli constants are by (\ref{q11}):

\newpage

\begin{align}
    \beta_1\, & =\, \pi \int_{0}^{\infty} \frac{a}{c^2} \, dz \,,  \label{q15}                                       \\
    \beta_2\, & =\,  \pi \int_{0}^{\infty} \frac{a^2 - b^2}{ c^2} \,dz \,,  \label{q16}                       \\
    \beta_3\, & =\, \pi \int_{0}^{\infty} \frac{a^3 - 3 a b^2 }{ c^2} \,dz \,, \label{q17}                  \\
    \beta_4\, & =\, \pi \int_{0}^{\infty} \frac{a^4 - 6 a^2b^2 + b^4}{c^2} \,dz \,, \label{q18}    \\
    \beta_5\, & =\, \pi \int_{0}^{\infty} \frac{a^5 - 10 a^3b^2 + 5ab^4}{c^2} \,dz \,. \label{q19}
\end{align}

\sect{Generalized Bernoulli constants}{20}

We recall that \( \gamma_n/n! \) is the coefficient of \( (1-s)^n \)
in the Laurent expansion of \( \zeta(s) \) about \( s = 1 \) and
\( \gamma_n(v)/n! \) is the coefficient of \( (1-s)^n \) in the Laurent
expansion of \( \zeta(s,v) \) about \( s = 1 \). In other words, with the
\textit{generalized Stieltjes constants}\index{Stieltjes!gen. constants} \( \gamma_n(v) \)
we have the \textit{Hurwitz zeta function}\index{Hurwitz!zeta function} \(\zeta(s,v)\) in the form
\begin{equation}
    \zeta(s,v)=\frac{1}{s-1}+\sum_{n=0}^\infty \frac{(-1)^n}{n!} \gamma_n(v) (s-1)^n, \quad(s \neq 1).
    \label{q35}
\end{equation}
The generalized Stieltjes constants may be computed for \( n \ge 0 \) and \( \Re(v) > \frac12 \) by an extension
of the integral representation (\ref{q2}), see Johansson and Blagouchine \cite[formulas~2, 32, 42]{Johansson}.
\begin{equation}
    \gamma_n(v) \, = \, - \frac{\pi}{2(n+1)}  \int_{- \infty}^\infty\frac{\log\!\left(v- \tfrac{1}{2} + i z\right) ^{n+1}}
    {\cosh(\pi x)^2}\, dz.
    \label{int333}
\end{equation}
The \textit{generalized Bernoulli constants}\index{Bernoulli!general constants} are defined as
\begin{equation}
    {
        \beta_s(v) = 2 \pi \!\int_{- \infty}^{+\infty}
        \!\frac{\, \log\!\big(v - \frac12 + ix\big) ^s \,}{ \left( {\Ee}^{- \pi x} + {\Ee}^{ \pi x} \right)^2 }\,dx.
    }
    \label{q37}
\end{equation}
We see that  \(\beta_{0}(v) = 1\) for all \(v\)  and
\( \beta_{n}(v) = - n \gamma_{n-1}(v) \) for \( n \ge 1 \).

\sect{The generalized Bernoulli function}{20}

\noindent We introduce the generalized Bernoulli function \( \operatorname{B}(s, v)\)
analogous to the Hurwitz zeta function. The new  parameter \( v \) can be
any complex number that is not a nonpositive integer. Then we define the
\textit{generalized Bernoulli function}\index{Bernoulli!generalized function} as
\begin{equation}{
        \operatorname{B}(s,v) \, = \, \sum _{n=0}^{\infty }\beta_{n}(v) \frac {s^n}{n!}\, .}
    \label{Generalizedgamma}
\end{equation}
For \( v=1, \) this is the ordinary Bernoulli function (\ref{Bdef}).
Using the identities from the last section, we get
\begin{equation}
    \operatorname{B}(s, v) = 1 - s \sum_{n=0}^{\infty} \gamma_n(v) \frac{s^n}{n!} .
    \label{q36a}
\end{equation}
Thus the generalized Bernoulli function can be represented by
\begin{equation}
    \operatorname{B}(s, v) = - s \, \zeta(1-s, v), \qquad (s \ne 1).
    \label{q36}
\end{equation}
This identity also embeds the \textit{Bernoulli polynomials}\index{Bernoulli!polynomials} as
\begin{equation}
    \operatorname{B}_{n}(x) \,=\, \operatorname{B}(n, x)  \quad (n \ge 0, n \text{ integer}) .
    \label{q22}
\end{equation}
This follows from (\ref{q36}) (see for instance Apostol \cite[th.\ 12.13]{Apostol})
and the fact that \(\operatorname{B}(0, x) \, = \, 1.\)

\sect{Integral formulas for the Bernoulli function}{20}

\noindent We can transfer the integral formulas for the Bernoulli constants to the
Bernoulli function itself. First we reproduce a formula by Jensen \cite{Jensen},
which he gave in reply to Cesàro in \textit{L'Interm\'ediaire des math\'ematiciens}.
\begin{equation} {
    (s-1)\zeta(s) \, = \, 2 \pi \int_{- \infty}^\infty\frac{( \tfrac{1}{2} +iy)^{1-s}}{({\Ee}^{\pi y} + {\Ee}^ {- \pi y})^2}\, dy.
    \label{int1}
    }\end{equation}
Jensen comments:
\begin{quote}
    \textit{``... [this formula] is remarkable because of its simplicity
    and can easily be demonstrated with the help of Cauchy's theorem.''}
\end{quote}
How Jensen actually computed \((s-1)\zeta(s)\) is unclear, since the formula for
the coefficients \(\operatorname{c}_v\), which he states, rapidly diverges.
This was observed by Kotěšovec (personal communication).

In a numerical example Jensen uses the Bernoulli constants in the form
\(\operatorname{c}_v \, = \, (-1)^v \beta_v/ {v!}. \)
Applied to the Bernoulli function, \textit{Jensen's formula}\index{Jensen!formula} is written as
\begin{equation} {
    \operatorname{B}(s) \, = \, 2 \pi \int_{- \infty}^\infty\frac{( \tfrac{1}{2} +iz)^{s}}
    {({\Ee}^{\pi z} + {\Ee}^ {- \pi z})^2}\, dz.
    \label{int4}
    } \end{equation}
In Johansson and Blagouchine \cite{Johansson} this is a particular case of the first
formula of theorem 1.  (See also Srivastava and Choi \cite[p.\ 92]{Srivastava}
and the discussion \cite{Stopple}.)

Hadjicostas \cite{Hadjicostas} remarks that from this theorem also the
representation for the generalized Bernoulli function follows:

\medskip
\textit{For all} \( s \in \mathbb{C}\) \textit{and} \( v \in \mathbb{C}\) \textit{with} \(\Re(v) \ge 1/2\)
\begin{equation} {
    \operatorname{B}(s, v) \, = \, 2 \pi \int_{- \infty}^\infty\frac{(v- \tfrac{1}{2} + i z)^{s}}
    {({\Ee}^{ \pi z} + {\Ee}^ {- \pi z})^2}\, dz.
    \label{int5}
    } \end{equation}
This formula is the central formula in this paper.

\begin{figure}
    \centering
    \includegraphics[height=185pt,width=310pt]{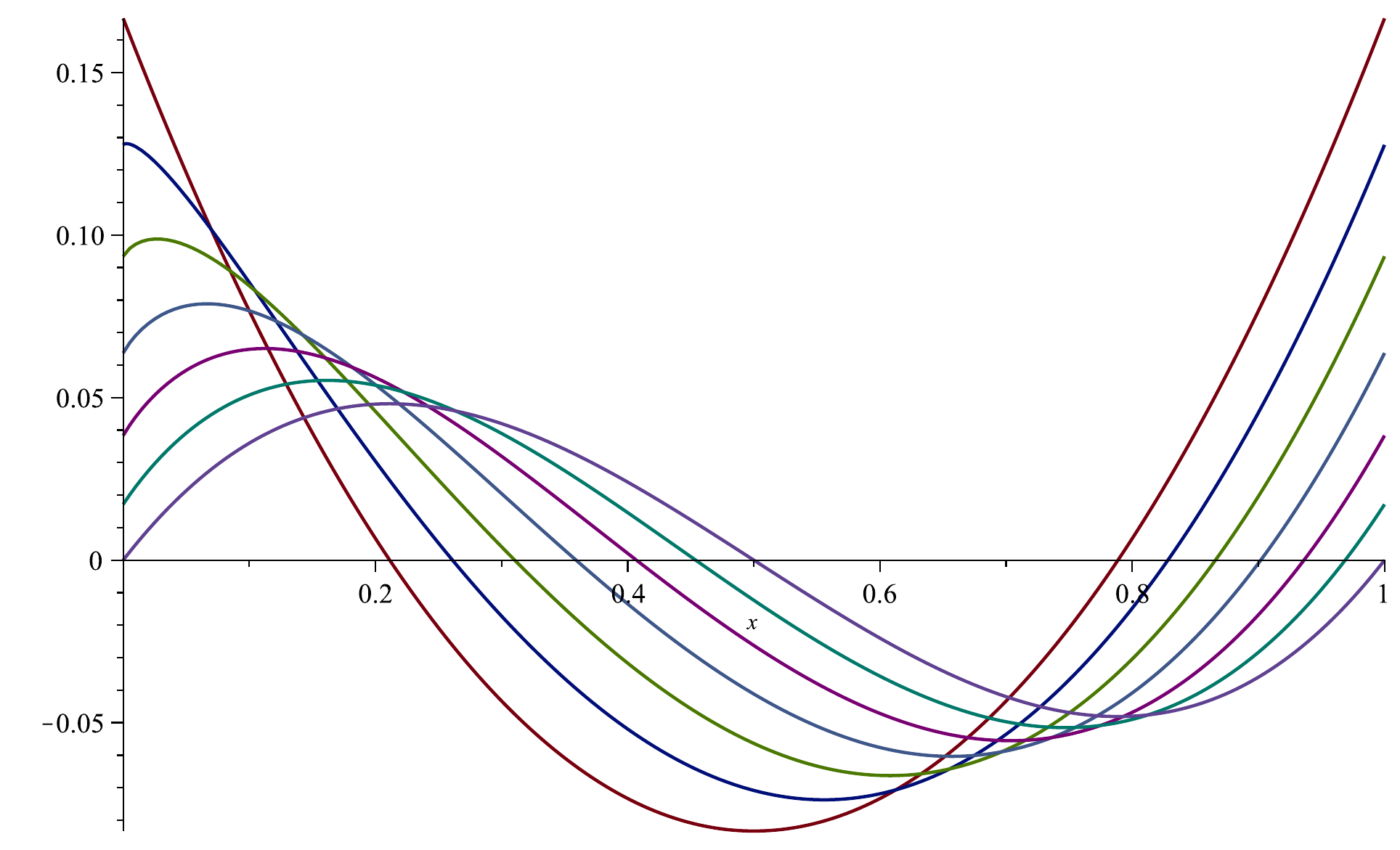}%
    \caption[The Hurwitz Bernoulli function]{\centerline{\textit{The Hurwitz Bernoulli functions with \( s = 2 + k/6,\) }}
  \centerline{\((0 \le k \le 6), \) \textit{deform \(\operatorname{B}_{2}(x)\) into \(\operatorname{B}_{3}(x).\)} }}%
  \label{fig:fig4hurbern}
\end{figure}

\sect{The Hurwitz--Bernoulli function}{20}

\noindent The \textit{Hurwitz--Bernoulli function}\index{Bernoulli!Hurwitz function} is defined as
    {
        \begin{align} \begin{split} \operatorname{H}(s, v) &= {\Ee}^{-i \pi s/2}\,
                \operatorname{L}(s, v) + {\Ee}^{i \pi s/2} \, \operatorname{L}(s, 1-v),   \\
                \operatorname{L}(s, v) &= - \frac{s\, !}{(2 \pi)^s} \operatorname{Li}_s({\Ee}^{2 \pi i v}). \label{hurbern}
            \end{split} \end{align}
   }
Here \( \operatorname{Li}_s(v) \) denotes the polylogarithm.
The proposition that
\begin{equation}
    \operatorname{B}_{s}(v) =\operatorname{B}(s, v) =
    \operatorname{H}(s, v), \  \  \text{for }  0 \le v \le 1  \text{ and }  s > 1,
    \label{HurwitzFormula}
\end{equation}
goes back to Hurwitz. The corresponding case for the zeta function (\ref{HurwitzFormula})
is known as the \textit{Hurwitz formula}\index{Hurwitz!formula} \cite[p.\ 71]{Apostol}.
With the Hurwitz--Bernoulli function the Bernoulli polynomials can be
continuously deformed into each other (see figure \underline{\ref{fig:fig4hurbern}}).

\begin{figure}
    \centering
    \includegraphics[width=310pt,height=180pt]{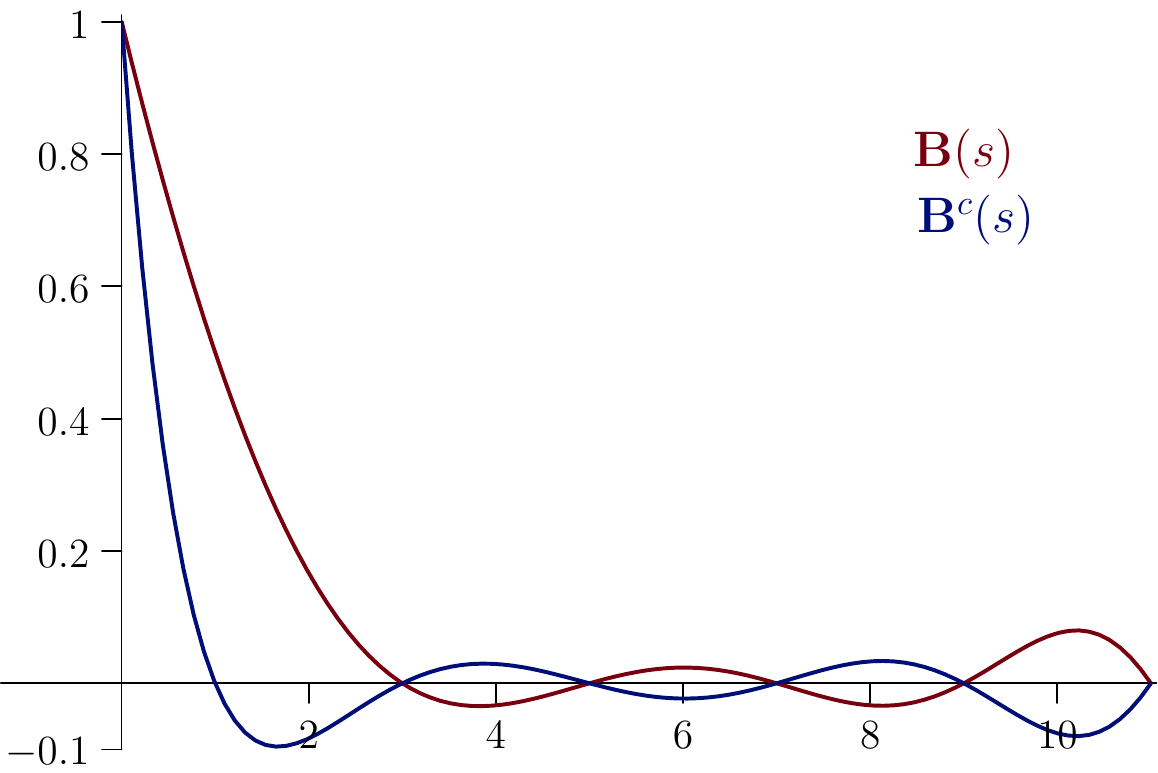}%
    \caption[The central Bernoulli function]{\textit{The Bernoulli function and the central Bernoulli function}}%
        \label{fig:fig5centbern}
\end{figure}

\sect{The central Bernoulli function}{20}

\noindent Setting \(v=1\) in (\ref{HurwitzFormula}) the Hurwitz--Bernoulli function simplifies to
\begin{equation} {
        \operatorname{B}_{s}(1) \, =\, -2\, s!  \operatorname{Li}_s(1) \cos(s \pi /2)/ (2 \pi)^{s}, \quad (s>1).
        \label{bs1}
    } \end{equation}
For \(s>1\) one can replace the polylogarithm with the zeta function and then apply
the functional equation of the zeta function to get
\begin{equation}
    \operatorname{B}_{s}(1) \, =\, -s \zeta(1-s),  \qquad (s>1).
    \label{brep}
\end{equation}
Thus the Bernoulli function is a vertical section of the Hurwitz--Bernoulli function,
\(\operatorname{B}(s)  = \operatorname{B}_{s}(1), \) similarly as the Bernoulli numbers
are special cases of the Bernoulli polynomials, \(\operatorname{B}_n = \operatorname{B}_{n}(1) \).

Setting \(v=1/2\) in the Hurwitz--Bernoulli function leads to a second noteworthy case.
Then (\ref{HurwitzFormula}) reduces to
\begin{equation} {
        \operatorname{B}_{s}(1/2) \, =\, -2\, s!
        \operatorname{Li}_s(-1) \cos(s \pi /2)/ (2 \pi)^{s}, \quad (s>1).
        \label{bs12}
    } \end{equation}
For \(s>1\) we can replace the polylogarithm with the negated alternating zeta function.
We call {\(\operatorname{B}^{c}(s)  = \operatorname{B}_{s}(1/2) \)} the
\textit{central Bernoulli function}\index{Bernoulli!central function}.

\newpage
The function can be expressed as \( \operatorname{B}^{c}(s)  =  -s\zeta(1-s)(2^{1-s} - 1) \) for \(s > 1 \), or
by the Bernoulli function as
\begin{equation}
    \operatorname{B}^{c}(s)  = \operatorname{B}(s)(2^{1-s} - 1).
    \label{cbrep}
\end{equation}

The central Bernoulli function has the same trivial zeros as the Bernoulli
function, plus a zero at the point \(s=1\) (see  figure \underline{\ref{fig:fig5centbern}}).

\smallskip
An integral representation for the central Bernoulli function follows from (\ref{int5}).
For all \( s \in \mathbb{C}\)
\begin{equation} {
        \operatorname{B}^{c}(s) \, = \, 2 \pi \int_{- \infty}^\infty\frac{(i z)^{s}}
        {({\Ee}^{- \pi z} + {\Ee}^ {\pi z})^2}\, dz.
        \label{intcenter5}
} \end{equation}
Assuming \(s >0\) and real, it follows from  (\ref{cbrep}) and (\ref{intcenter5}) 
that%
\begin{equation} {
        \operatorname{B}(s) \, = \, \pi \,\frac{\cos(\pi s / 2)}{2^{1 - s} - 1}
        \int_{0}^\infty z^{s} \sech(\pi z)^2 \, dz.
        \label{intrealbern}
} \end{equation}

For \(s = 1\) the value on the right-hand side of
(\ref{intrealbern}) is to be understood as the limit value \((\pi^2/\log(4)) (\log(2)/\pi^2) = 1/2\).

\smallskip
We call  (\ref{intrealbern}) the \textit{secant decomposition} of  \(\operatorname{B}(s)\) (see figure \underline{\ref{fig:fig35SecDecomp}}).

\begin{figure}
    \centering
    \includegraphics[height=240pt,width=320pt]{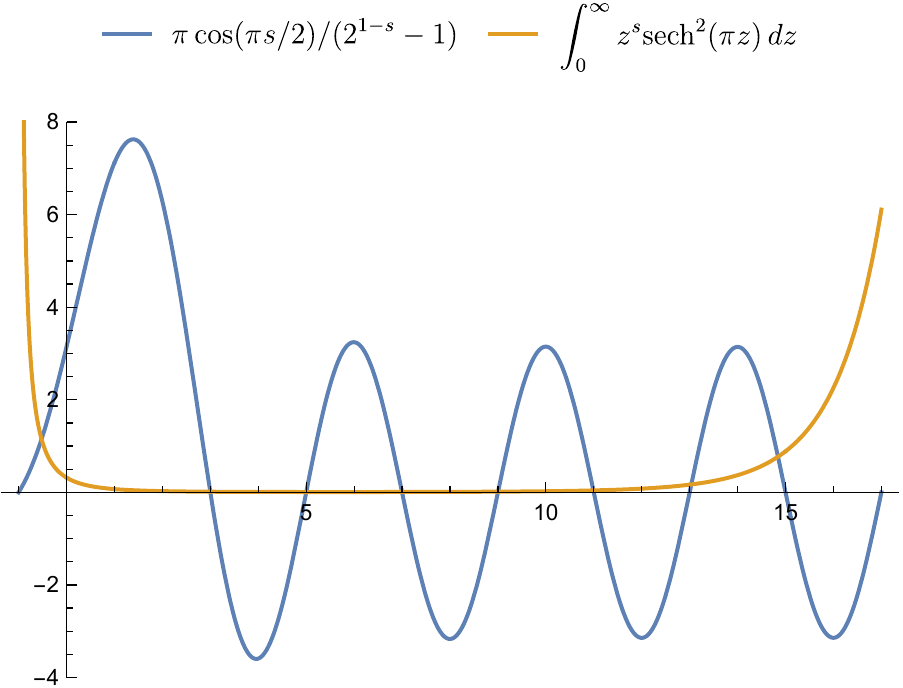}%
    \medskip
    \caption[The secant decomposition]{\textit{The secant decomposition of the real Bernoulli function.}}
    \label{fig:fig35SecDecomp}
\end{figure}    

\FloatBarrier
\newpage

\sect{The central Bernoulli polynomials}{0}

\noindent The \textit{central Bernoulli numbers} are defined as
\begin{equation}
    \operatorname{B}^{c}_n  = 2^n \operatorname{B}_{n}(1/2).
    \label{CBN}
\end{equation}
The first few are, for \(n \ge 0\):
\begin{equation}
    1,\, 0,\, - \frac{1}{3},\, 0,\, \frac{7}{15},\, 0,\, - \frac{31}{21},\, 0,\, \frac{127}{15},
    \, 0,\, - \frac{2555}{33}, \, 0, \frac{1414477}{1365}, \ldots .
    \nonumber
\end{equation}

\begin{table}[t]
    \centering \setlength{\extrarowheight}{3pt}
    {  \begin{tabular}{c}
            \( \phantom{-}1 \)                                                                                                                   \\
            \( \phantom{-} x \)                                                                                                                   \\
            \( -{\frac{1}{3}}\,+\,{x}^{2} \)                                                                                                  \\
            \( -{x}\,+\,{x}^{3}\)                                                                                                                    \\
            \( \phantom{-}{\frac{7}{15}}\,- \,{2\, {x}^{2}}\,+\,{x}^{4}\)                                               \\
            \( \phantom{-}{\frac {7}{3}\,x}\,- \,{\frac {10}{3}\,{x}^{3}}\,+\,{x}^{5}\)                         \\
            \( -{\frac{31}{21}}\,+\,{7\,{x}^{2}}\,- \,5\,{x}^{4}\,+\,{x}^{6}\)                                           \\
            \( -{\frac {31}{3}\,x}\,+\,{\frac {49}{3}\,{x}^{3}}\,- \,7\,{x}^{5}\,+\,{x}^{7}\)                    \\
            \( {\frac{127}{15}}-{\frac {124}{3}\,{x}^{2}}+{\frac {98}{3}\,{x}^{4}}-{\frac {28}{3}\,{x}^{6}}+{x}^{8}\) \\
            \( {\frac {381}{5}\,x}-124\,{x}^{3}+{\frac {294}{5}\,{x}^{5}}-12\,{x}^{7}+{x}^{9}\)
        \end{tabular} }
    \caption{\textit{The central Bernoulli polynomials} \(\operatorname{B}_n^{c}(x)\).}
    \label{centralbernpoly}
\end{table}

Unsurprisingly Euler, in 1755 in his \textit{Institutiones}, also calculated some
central Bernoulli numbers,  B(3,1) and B(5,1) (Opera Omnia, Ser. 1, Vol. 10, p.~351).

The \textit{central Bernoulli polynomials}\index{Bernoulli!central polynomials}\index{Polynomials!central Bernoulli}
are by definition the Appell sequence
\begin{equation}
    \operatorname{B}_n^{c}(x) \, = \,  \sum_{k=0}^n \binom{n}{k} \operatorname{B}^{c}_k x^{n-k}.
\end{equation}
The parity of \(n\) equals the parity of \(\operatorname{B}_n^{c}(x),\)
a property the Bernoulli polynomials do not possess.

Despite their systematic significance, the central Bernoulli polynomials were not in the
\OEIS database at the time of writing these lines (now they are filed in \OEIS \seqnum{A335953}).

Also the following identity is worth noting:
\begin{equation}
   2^n \operatorname{B}_n(1) \, = \,  \sum_{k=0}^n \binom{n}{k}  2^{k}
    \operatorname{B}_k\left( \tfrac{1}{2} \right).
    \label{bernoulligoeshalf}
\end{equation}
We will come back to this identity later, as it highlights the correct choice
of the generating function of the Bernoulli numbers.

 \begin{table}[t]
    \centering \setlength{\extrarowheight}{3pt}
    {  \begin{tabular}{c}
            \( 0 \) \\
            \(                 -1 \) \\
            \(                 1 - 2 x  \) \\
            \(                 3 x - 3 x^2 \) \\
            \(                 -1  +6 x^2 -4 x^3 \) \\
            \(                -5 x +10 x^3 -5 x^4 \) \\
            \(         3  -15 x^2   +15 x^4   -6 x^5 \) \\
            \(            21 x  -35 x^3  +21 x^5  -7 x^6  \) \\
            \(          -17   +84 x^2 -70 x^4 +28 x^6  -8 x^7 \)
    \end{tabular} }
    \caption{\textit{The Genocchi polynomials} \(\operatorname{G}_n(x)\).}
    \label{genocchipoly}
\end{table}%

\sect{The Genocchi function}{20}

How much does the central Bernoulli function deviate from the Bernoulli function?
The \textit{Genocchi function}\index{Genocchi!function} answers this question (up to a normalization factor).
\begin{equation}
    \operatorname{G}(s)  \,
    = \, 2^s \left( \operatorname{B}^{c}(s) - \operatorname{B}(s) \right)
      = \, 2^s \left( \operatorname{B}_{s}( \tfrac{1}{2} ) - \operatorname{B}_{s}(1) \right) .
    \label{GenocciFun}
\end{equation}
From the identities (\ref{brep}) and (\ref{cbrep}), it follows that the Genocchi function
can be represented as
\begin{equation}
    \operatorname{G}(s)  \,  = \,  2 (1 - 2^s) \operatorname{B}(s).
    \label{GenocciFun2}
\end{equation}
The Genocchi function takes integer values for nonnegative integer arguments.
The \textit{Genocchi numbers}\index{Genocchi!numbers}
 \( \operatorname{G}_n = \operatorname{G}(n) \) are  listed in \OEIS \seqnum{A226158}.
 The correct sign of \( \operatorname{G}_1 = -1 \) must be observed.
  \[
     \operatorname{G}_n  =   0,\, -1,\, -1,\, 0,\, 1,\, 0,\, -3,\, 0,\, 17,\, 0,\, -155,\, 0,\, 2073,\, \dots
 \]
The \textit{generalized Genocchi function}\index{Genocchi!generalized function} generalizes the definition (\ref{GenocciFun}),
\begin{equation}
    \operatorname{G}(s, v) \, = \, 2^s \left(\operatorname{B}\left(s, \frac{v}{2}\right) -
    \operatorname{B}\left(s, \frac{v+1}{2} \right) \right) \label{GG}.
\end{equation}
The plain Genocchi function is recovered by  \( \operatorname{G}(s) =  \operatorname{G}(s, 1) \).

\noindent For \(n\) integer (\(\ref{GG}\)) gives the
\textit{Genocchi polynomials}\index{Genocchi!polynomials}\index{Polynomials!Genocchi}
\begin{equation} {
        \operatorname{G}_{n}(x)
        = \, 2^n \left( \operatorname{B}_n\left(\frac{x}{2}\right) -
        \operatorname{B}_n\left( \frac{x+1}{2}\right) \right).}
    \label{q123}
\end{equation}
The Genocchi polynomials are twice  the difference between the Bernoulli polynomials
and the central Bernoulli polynomials.
\begin{equation} {
        \operatorname{G}_{n}(x)
        = \, 2 \left( \operatorname{B}_n(x)  -
        \operatorname{B}^{c}_n(x) \right).}
    \label{q3123}
\end{equation}
The integer coefficients of these polynomials (with different signs) are
 \seqnum{A333303} in the \OEIS. The Genocchi function is closely related to
the alternating Bernoulli function, as we will see next.

\begin{figure}
    \centering
    \includegraphics[height=240pt,width=320pt]{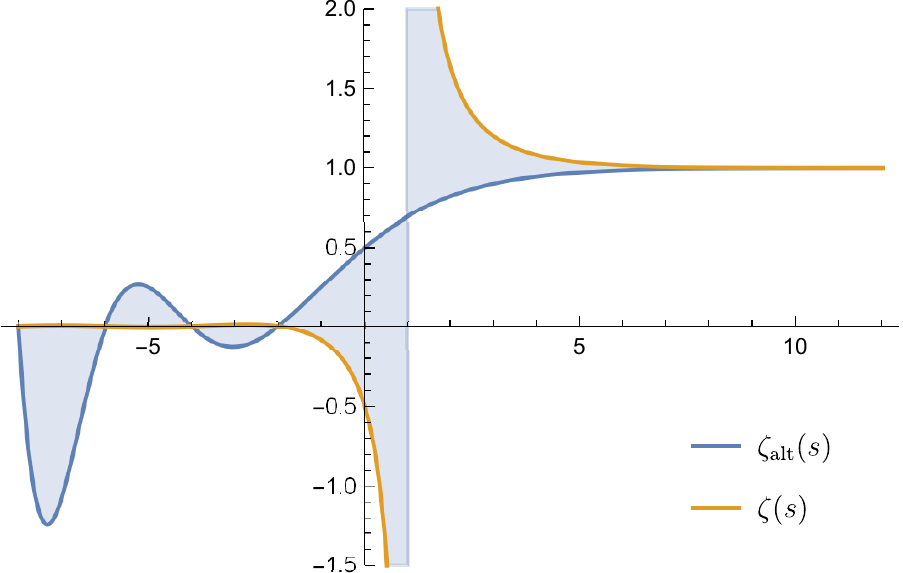}%
    \caption[Riemann alternating]{\textit{Riemann zeta versus Riemann alternating zeta.}} 
    \label{fig:fig25riemannalternating}
\end{figure}

\sect{The alternating Bernoulli function}{20}

\noindent The \textit{alternating Riemann zeta function}\index{alternating zeta!Riemann function},
also known as the \textit{Dirichlet eta function}\index{Dirichlet eta function}, is defined as
\begin{equation}
    \zeta^{*}(s) = (1- 2^{1-s})\zeta(s), \quad (s \ne 1).
    \label{q57}
\end{equation}
The \textit{alternating Bernoulli function}\index{alternating Bernoulli!function} is defined
in analogy to the representation of the Bernoulli function by the zeta function.
\begin{equation}
    \operatorname{B}^{*}(s) = -s \zeta^{*}(1-s), \qquad (s \ne 0).
    \label{betastardef}
\end{equation}
The \textit{alternating Hurwitz zeta function}\index{alternating zeta!Hurwitz function}
is set for \(0 < x \le 1\) as
\begin{equation}
    \zeta^{*}(s, x) \, =\, \sum_{n=0}^{\infty} \frac{(-1)^n}{(n+x)^s},  \quad \text{for } \Re(s) > 0,
    \label{a61}
\end{equation}
and for other values by analytic continuation. This function is represented by the Hurwitz Zeta
function as
\begin{equation}
    \zeta^{*}(s, x) \, =\, 2^{-s} \left(\zeta\left(s, \frac{x}{2}\right) - \zeta\left(s, \frac{x+1}{2}\right)  \right), \quad (s \neq 1).
    \label{b61}
\end{equation}
The \textit{alternating Bernoulli polynomials}\index{alternating Bernoulli! polynomials}
are defined as
\begin{equation}
    \operatorname{B}^{*}_n(x) \, = \, -n \zeta^{*}(1 - n, x)  \,.
    \label{bstar}
\end{equation}
This definition is in analogy to the introduction of the Bernoulli polynomials (\ref{q22}).

The alternating Bernoulli function can be represented by the Bernoulli function
\begin{equation}
    \operatorname{B}^{*}(s) =  \left( 1-2^{s} \right) \operatorname{B}(s).
    \label{betastar}
\end{equation}

\begin{figure}
    \centering
    \includegraphics[height=240pt,width=320pt]{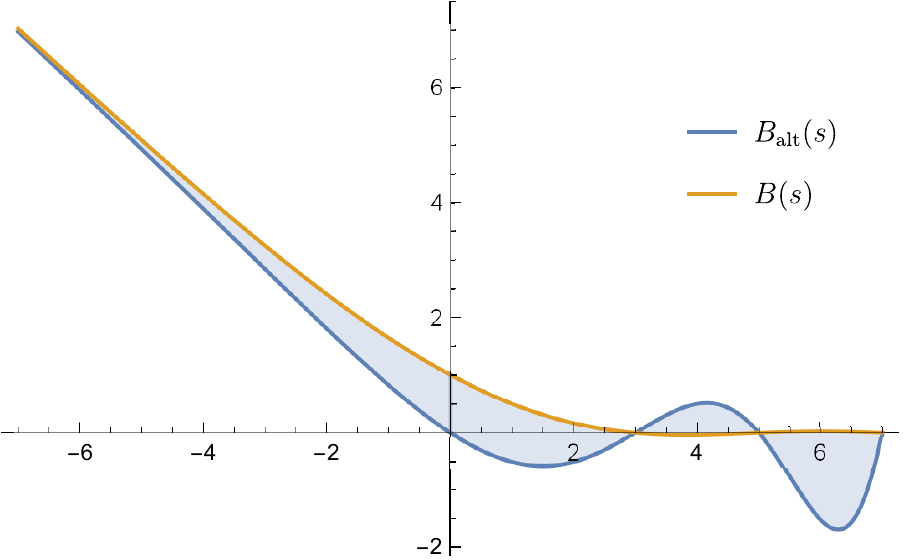}%
    \caption[Bernoulli alternating]{\textit{Bernoulli function versus Bernoulli alternating function.}} 
    \label{fig:fig26bernoullialternating}
\end{figure}
 \FloatBarrier

The \textit{alternating Bernoulli numbers}\index{alternating Bernoulli!numbers}
  \(  \operatorname{B}^{*}_n = \operatorname{B}^{*}(n)  \label{q61} \)
are the values of the alternating Bernoulli function at the nonnegative integers,
and are, like the Bernoulli numbers,  rational numbers.

Reduced to lowest terms, they have the denominator \( 2 \)
and by (\ref{q123}) are half the  Genocchi numbers.
\begin{equation}
    \operatorname{B}^{*}_{n}
    \, = \, 2^{n-1} \left( \operatorname{B}_n\left(\frac{1}{2}\right) - \operatorname{B}_n\left(1 \right) \right) =  \operatorname{G}_n/2.
\end{equation}

\begin{figure}
    \centering
    \includegraphics[width=340pt,height=205pt]{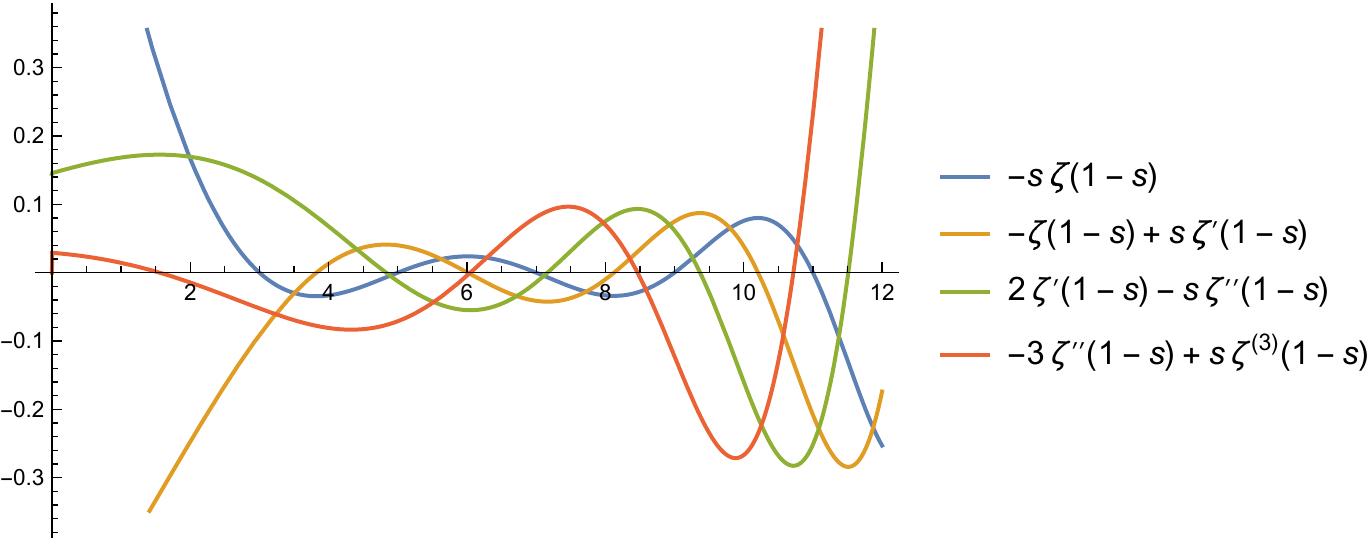}%
    \caption[Derivatives of the Bernoulli function]{\textit{The derivatives of the Bernoulli function
            \( \operatorname{B}^{(n)}(s),\, 0 \le n \le 3. \)}}
    \label{fig:fig7bernfuneulgam}    
\end{figure}

\sect{Derivatives of the Bernoulli function}{20}

\noindent The Bernoulli constants are in a simple relationship with the derivatives of the Bernoulli function.
\index{Bernoulli!derivative} With the Riemann zeta function we have
\begin{equation} {
        \operatorname{B}^{(n)}(s) = (-1)^n \left( n \zeta^{(n-1)}(1-s)
        - s \zeta^{(n)}(1-s) \right). }
    \label{q25}
\end{equation}
Here \(\operatorname{B}^{(n)}(s) \) denotes the \(n\)-th derivative of the Bernoulli function.

Taking \( \lim_{s \rightarrow 0} \) on the right hand side of (\ref{q25}) we get
\begin{equation}
    \operatorname{B}^{(n)}(0) = \beta_n  = - n \, \gamma_{n-1}  \quad ( n \ge 1) \, .
    \label{q26}
\end{equation}

Entire books \cite{Havil} have been written about the emergence of
Euler's gamma in number theory. The identity\index{Constant!\(\gamma\ \)Euler}
\({- \operatorname{B}'(0) = \gamma }\label{q39} \)
is one of the beautiful places where this manifests (see figure \underline{\ref{fig:fig7bernfuneulgam}}).

\begin{figure}
    \centering
    \includegraphics[width=320pt]{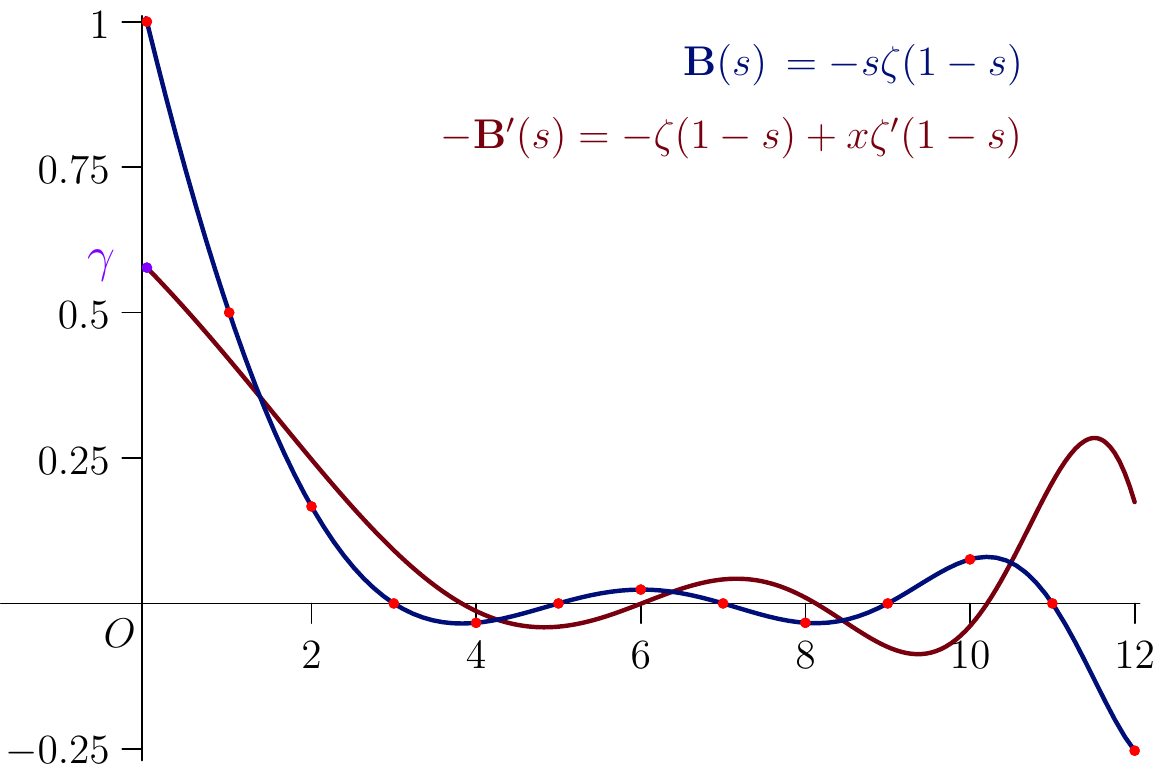}%
    \caption[Euler's \(\gamma\) and \(\operatorname{B}'(x)\)]{\textit{\( -\operatorname{B}'(x)\) hits Euler's \( \gamma\) at  \( x=0, \) dots are Bernoulli numbers.}}
\end{figure}

\sect{Logarithmic derivative and  Bernoulli cumulants}{0}

\noindent The \textit{logarithmic derivative of a function}
\index{logarithmic!derivative} \(\operatorname{F}\) will be denoted by
\begin{equation}
    \mathcal{L}\!\operatorname{F}(s) = \frac{ \operatorname{F}'(s)}{ \operatorname{F}(s)}.
    \nonumber
\end{equation}

In particular we will write \(\mathcal{L}\!\operatorname{B}(s),\, \mathcal{L}\zeta(s)\), and
\(\mathcal{L}\Gamma(s)\) for the logarithmic derivative of the Bernoulli function,
\index{Bernoulli!log. derivative} the \(\zeta\) function, and the \(\Gamma\) function.
\(\mathcal{L}\Gamma(s)\) is also known as the digamma function \(\psi\).

In terms of the zeta function \(\mathcal{L}\!\operatorname{B}(s)\) can also be written
\begin{equation}
    \mathcal{L}\!\operatorname{B}(s)
    \,=\, \frac{1}{s} - \mathcal{L}\zeta(1-s).
    \label{ws2}
\end{equation}
If \(s=0\) then \(\mathcal{L}\!\operatorname{B}(0) \) is set to the limiting value  \(-\gamma \).
For odd integer \(n = 3, 5, 7, \ldots \) the value of  \(\mathcal{L}\!\operatorname{B}(n) \) is undefined.

From (\ref{ws2}) we can infer
\begin{equation} {
        \mathcal{L}\!\operatorname{B}(s)
        \,=\,   \mathcal{L}\Gamma(s)  + \mathcal{L}\zeta(s) + \rho(s) },
    \label{BZ}
\end{equation}
where  \( \rho(s) \, =\,  1/s  - (\pi /2) \tan(s \pi / 2) - \log(2\pi). \)

\medskip
The series expansion  of
\( \mathcal{L}\!\operatorname{B}(s)  =  \sum_{n = 1}^{\infty} b_n s^{n-1} \)  at \( s = 0 \)  starts
\begin{align*}
     & \mathcal{L}\!\operatorname{B}(s) = \beta_{{1}}
    +( \beta_{{2}} - {\beta_{{1}}}^{2}  ) \,s
    +( \beta_{{3}} - 3\,\beta_{{1}}\beta_{{2}} + 2\,{\beta_{{1}}}^{3} ) \,s^2/2                                                                       \\
     & + ( \beta_{{4}} - 3\, {\beta_{{2}}}^{2} - 4\,\beta_{{1}}\beta_{{3}}
     +12\,{\beta_{{1}}}^{2}\beta_{{2}} -6\,{\beta_{{1}}}^{4} ) \,s^3/6  + O(s^4).
    \label{q41}
\end{align*}
The coefficients \(b_n \) are given for \( n \ge 1 \) by
\begin{equation}
    b_n \, = \, n! \,[ s^n ]   \log \bigg( \sum_{n = 0}^{\infty} \beta_n  \frac{s^n}{n!} \bigg).
    \label{q42}
\end{equation}
In other words, the coefficients are the \textit{logarithmic polynomials}\index{logarithmic!polynomials}
generated by the Bernoulli constants, Comtet \cite[p.\ 140]{Comtet}.
These polynomials may be called \textit{Bernoulli cumulants}\index{Bernoulli!cumulants},
following a similar naming by Voros \cite[3.16]{Voros}.
The numerical values appearing in this expansion,
listed as an irregular triangle, are \OEIS \seqnum{A263634}.

\bigskip

\begin{table}[hb]    \centering
    \begin{tabular}{cccccc}
        1  &                &             &            &     &                     \\
        -1 & 1            &             &            &     &                     \\
        2  & -3            & 1           &            &     &                    \\
        -6 &  12         & - 4,\  -3           & 1         &     &                \\
        24 & -60      & 20, \ 30      & -5,\ -10     & 1  &         \\
        -120 & 360 & -120, \ -270 & 30, \ 120, \ 30 & -6, \ -15, \ -10 & 1
    \end{tabular}
    \caption{\textit{Coefficients of the Bernoulli cumulants.}}
    \label{logpoly}
\end{table}

\sect{The Hasse--Worpitzky representation}{20}

\noindent The coefficients of the Bernoulli cumulants  are refinements
of the signed \textit{Worpitzky numbers}\index{Worpitzky!numbers}
\( \operatorname{W}(n,k) \) \cite{Worpitzky, Vandervelde}. See 
figure \underline{\ref{fig:fig30fubinipoly}} for the 
\textit{Worpitzky and Fubini polynomials,}\index{Worpitzky!polynomials}
\index{Fubini!polynomials} \OEIS \seqnum{A163626} and \seqnum{A278075}.
\begin{equation}
    \operatorname{W}(n,k) \,= \, (-1)^k  k! \stirling{n+1}{k+1}.
    \label{q43}
\end{equation}
Here \( \stirling{n}{k} \) denotes the Stirling set numbers.
Generalizations based on Joffe's central differences
of zero are \seqnum{A318259} and \seqnum{A318260}.

The \textit{Worpitzky transform}\index{Worpitzky!transform} maps a sequence
\(a_0, a_1, a_2, \ldots \) to a sequence \( b_0, b_1, b_2, \ldots, \)
\begin{equation}
    b_n \,= \, \sum_{k=0}^n \,  \operatorname{W}(n,k)\, a_k \,.
    \label{q44}
\end{equation}
If \(a\) has the ordinary generating function \(a(x)\), then \(b\) has
exponential generating function \(a(1- \mathrm{e}^x) \mathrm{e}^x \).
Merlini et al. \cite{Merlini}  call the transform the Akiyama--Tanigawa
transformation; in the \OEIS also the term
Bernoulli--Stirling transform is used.

Worpitzky proved in 1883 that if we choose \(a_k=\frac{1}{k+1}\) and
apply the transform (\ref{q44}), the result is the sequence of the Bernoulli numbers. This approach can be generalized.

\enlargethispage{12pt}
\newpage

\begin{figure}
    \centering
    \includegraphics[height=240pt,width=324pt]{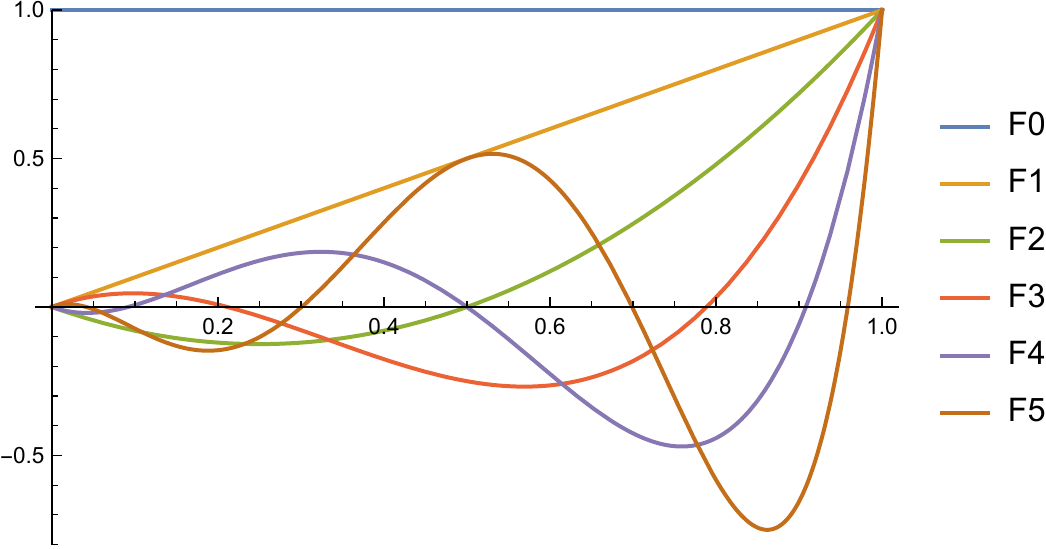}%
    \medskip
    \caption[Fubini polynomials]{\textit{Fubini and Worpitzky polynomials.}} 
    \index{Worpitzky!polynomials} \index{Fubini!polynomials} \index{Polynomials!Worpitzky}
    \vspace{-6pt}
    \caption[]{ \(\operatorname{F}_n(x)  = \sum_{k=0}^n (-1)^{n-k} k! \stirling{n}{k} x^k. \)} 
    \vspace{-6pt}
    \caption[]{\( \int_{0}^{1} \operatorname{F}_n(x) \, dx  = \operatorname{B_n} =  \int_{0}^{1} \operatorname{F}_n(1-x) \, dx.\)}
    \label{fig:fig30fubinipoly}
\end{figure}

\sect{The generalized Worpitzky transform}{10}

\noindent The \textit{generalized Worpitzky transform}\index{Worpitzky!gen. transform} 
\(\operatorname{W} : \mathbb{Z}^{\mathbb{N}} \rightarrow
\mathbb{Z}[x]^{\mathbb{N}}, \  a \in \mathbb{Z}^{\mathbb{N}}, \)
maps an integer sequence \(a_0, a_1, a_2, \ldots \)
to a sequence of polynomials.
The \(m \)-th term of  \(\operatorname{W}(a)\)  is the polynomial
\( \operatorname{W}_m(a)\), given by
\begin{equation}{
    \operatorname{W}\!\!_{m}(a) \, = \,
    \sum_{n=0}^m (-1)^n \binom{m}{n} x^{m-n}
    \sum_{k=0}^n \operatorname{W}(n,k) a_k .}
    \label{q45}
\end{equation}
In (\ref{q45}) the inner sum is the Worpitzky transform (\ref{q44}) of \(a\).
The first few polynomials\index{Polynomials!Worpitzky},
\( \operatorname{W}_0,  \ldots, \operatorname{W}_4, \) are:

\noindent  { \small
    \!\begin{minipage}[left]{320pt}
        \mathleft
        \begin{multline*}
            \mathtt{a_{{0}}}; \\
            \mathtt{a_{{0}}x - (a_{{0}} - a_{{1}}}); \\
            \mathtt{a_{{0}}{x}^{2} -  2\!\left( a_{{0}} - a_{{1}} \right)  x         +          \left(a_{{0}} - 3a_{{1}}+2a_{{2}} \right) };  \\
            \mathtt{a_{{0}}{x}^{3} -  3\!\left( a_{{0}} - a_{{1}} \right)  {x}^{2} +  3\! \left( a_{{0}} - 3a_{{1}} + 2a_{{2}} \right) x - \left(a_{{0}} - 7a_{{1}} + 12a_{{2}} - 6a_{{3}} \right) }; \\
            \mathtt{a_{{0}}{x}^{4} - 4\!\left( a_{{0}} - a_{{1}} \right)  {x}^{3} +  6\! \left( a_{{0}} - 3a_{{1}} + 2a_{{2}} \right) {x}^{2} - 4\! \left( a_{{0}}- 7a_{{1}} + 12a_{{2}} - 6a_{{3}}  \right)\! x } \\
            \mathtt{ + \left( a_{{0}}-15a_{{1}}+50 a_{{2}}-60a_{{3}}+24a_{{4}} \right). }
        \end{multline*}
    \end{minipage}
}

\smallskip
\noindent The definition (\ref{q45}) can be rewritten as
\begin{equation}
    \operatorname{W}\!\!_{m}(a) \, = \, \sum_{n=0}^m a_n \sum_{k=0}^n (-1)^{k} \binom{n}{k}  (x-k-1)^m .
    \label{q46}
\end{equation}

As the reader probably anticipated, we get the Bernoulli polynomials if we set
\(a_n = 1/(n+1)\). Evaluating at \(x=1\) we arrive at a well-known
representation of the Bernoulli numbers:
\begin{equation}
    \operatorname{B}_m \, = \,  \sum_{n=0}^m \frac{1}{n+1}
    \sum_{k=0}^n (-1)^{m-k} \binom{n}{k} k^m, \quad m \ge 0.
    \label{q47}
\end{equation}

It might be noted that if we choose
  \(a_n = n+1\) in  (\ref{q46}) then the  binomial polynomials are obtained.
The reader may also enjoy setting \( a_n = {H_{n+1}} \), where \({H}_n\) denotes the harmonic numbers.

\sect{The Hasse representation}{20}

\noindent In 1930 Hasse \cite{Hasse}  took the next step
in the development of formula (\ref{q47}) and proved:
\begin{equation}{
        \operatorname{B}(s,v)\,
        =\,\sum_{n=0}^\infty\frac{1}{\,n+1\,}
        \sum_{k=0}^n (-1)^k \binom{n}{k}(k+v)^{s}}
    \label{Hasse}
\end{equation}
 \textit{is an infinite series}\index{Hasse representation}
\textit{that converges for all complex \( s \)  and represents the entire
    function \( - s\zeta(1-s,v) \), the Bernoulli function.}

\smallskip
\noindent For instance the Hasse series of the central
Bernoulli function \index{Polynomials!central Bernoulli} is
\begin{equation}{
        \operatorname{B}^{c}(s) \,
        =\,  \sum_{n=0}^\infty\frac{1}{\,n+1\,}
        \sum_{k=0}^n (-1)^k \binom{n}{k} \left(k + \frac12 \right)^{s}}.
    \label{HasseCentral}
\end{equation}
This representation in turn provides an explicit formula for the central Bernoulli numbers
\index{central Bernoulli numbers}
(\( \ref{CBN} \)),
\begin{equation}{
        \operatorname{B}^{c}_n \,
        =\,  \sum_{j=0}^n \frac{1}{\,j+1\,}
        \sum_{k=0}^j (-1)^k \binom{j}{k} (2 k + 1)^{n}}.
    \label{HasseCentralNumbers}
\end{equation}
An important corollary to Hasse's formula is: The generalized Bernoulli constants \(\beta_s(v) \) \index{Bernoulli!general constants} (\ref{q37})
are given by the infinite series
\begin{equation}
    \beta_{s}(v) \,=\, \sum_{n=0}^\infty\frac{1}{\,n+1\,
    }\sum_{k=0}^n (-1)^k \binom{n}{k} \ln(k+v)^{s} .
    \label{q49}
\end{equation}
For proofs see Blagouchine \cite[Cor.~1, formula 123]{BlagouchineHasse} and the references given there.

\enlargethispage{12pt}
\newpage

\sect{The functional equation}{10}

\noindent The functional equation of the Bernoulli function
generalizes a formula of Euler which Graham et al. \cite[eq.\ 6.89]{ConcreteMath}
call \textit{almost miraculous}. But before we discuss it we introduce yet another
function and quote a remark from Tao \cite{TaoTau}.
\begin{quote}
    \textit{``It may be that \( 2 \pi i \) is an even more fundamental constant than
        \( 2 \pi \) or \( \pi \). It is, after all, the generator of \( \log(1) \). The fact
        that so many formulas involving \( \pi^n \) depend on the parity of \( n \)
        is another clue in this regard.''}
\end{quote}
Taking up this remark we will use the notation \( \tau  = 2 \pi i \)\index{Constant!\(\tau\)} and
the function \index{Tau function}
\begin{equation}
    {\tau(s)}  = \tau^{-s} + ( - \tau)^{-s} .
    \label{q50}
\end{equation}
We use the principal branch of the logarithm when taking powers of \(\tau\).
If \(s\) is real we can also write \(\tau(s)= 2^{1 - s} \pi^{-s} \cos(\pi s/2)\).%

This notation allows us to express the functional equation of
the Riemann zeta function \cite{RiemannEng} as the product of three functions,
\begin{equation}
    \zeta\!\left(1 - s\right) = \zeta\!\left(s\right)  \tau(s) \,  \Gamma(s)
    \qquad (s \in \mathbb{C} \setminus \{0, 1 \} ).
    \label{q51}
\end{equation}
Using (\ref{q51}) and \( \operatorname{B}(s)=-s\zeta(1-s) \) we get the representation
\begin{equation}{
        \operatorname{B}(s) = -  \zeta(s)\, \tau(s)\, s! .}
    \label{q52}
\end{equation}
From \( \operatorname{B}(1-s) = (s-1) \zeta(s) \) we obtain a self-referential
representation of the Bernoulli function, the \textit{functional equation}\index{Bernoulli!functional equation}
\begin{equation} {
        \operatorname{B}(s) =  \frac{\operatorname{B}(1-s)}{1-s}
        \tau(s)\, s! \, .}
    \label{q53}
\end{equation}
This functional equation also has a \textit{symmetric variant},
which means that the left side of (\ref{q54}) is unchanged by
the substitution \( s \leftarrow 1-s  \).%
\begin{equation}
    \operatorname{B}(1-s) \, \left( \frac{s}{2} \right) ! \ \pi^{-s/2} \,  = \,
    \operatorname{B}(s) \left( \frac{1-s}{2} \right)! \ \pi^{-(1-s)/2}.
    \label{q54}
\end{equation}

\newpage

\begin{figure}
    \includegraphics[height=240pt,width=320pt]{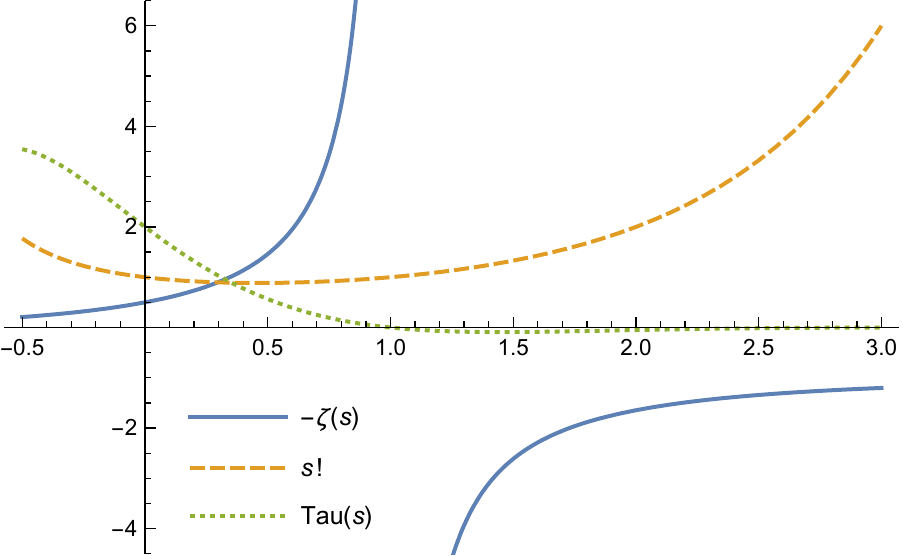}%
    \medskip
    \caption[The Riemann decomposition]{\textit{The Riemann decomposition of the Bernoulli function (\ref{q52}).}}
    \label{fig:fig36RiemannDecomp}
\end{figure}

\sect{Representation by the Riemann \(\xi\) function}{0}

\noindent The right (or the left) side of (\ref{q54}) turns out to be the
\textit{Riemann \(\xi\) function}\index{Riemann!\(\xi\) function},
\begin{equation}
    \xi(s) \,  = \, \left( \frac{s}{2} \right) !  \ \pi^{-s/2}\, (s-1) \zeta(s).
    \label{t35}
\end{equation}
For a discussion of this function see for instance Edwards \cite{Edwards}.
Thus we get a second representation of the Bernoulli function
in terms of a Riemann function:
\begin{equation}
    { \operatorname{B}(s)  \,  = \,
        \frac{ \pi^{(1-s)/2}}{\left( (1-s)/2 \right)! } \ \xi(s) \,. } \label{q554}
\end{equation}
By the functional equation of \( \xi \), \( \xi(s) = \xi(1-s) \), we also get
\begin{equation}
    { \operatorname{B}(1-s)  \,  = \,
        \frac{ \pi^{s/2} }{ \left( s/2 \right)! } \ \xi(s) \,. } \label{q557}
\end{equation}
This is a good opportunity to check the value of \(\operatorname{B}(-1)\).
\begin{equation}
    { \operatorname{B}(-1)  \,  = \,  \pi \,  \xi(-1) = \,\frac{ \pi^2}{6 }  = \pi \, \xi(2) \,.}
    \label{wq5}
\end{equation}
In this row of identities, the names Bernoulli, Euler (solving the Basel problem in 1734),
and Riemann join together.

\enlargethispage{12pt}
\newpage

\begin{figure}
    \centering
    \includegraphics[height=190pt,width=310pt]{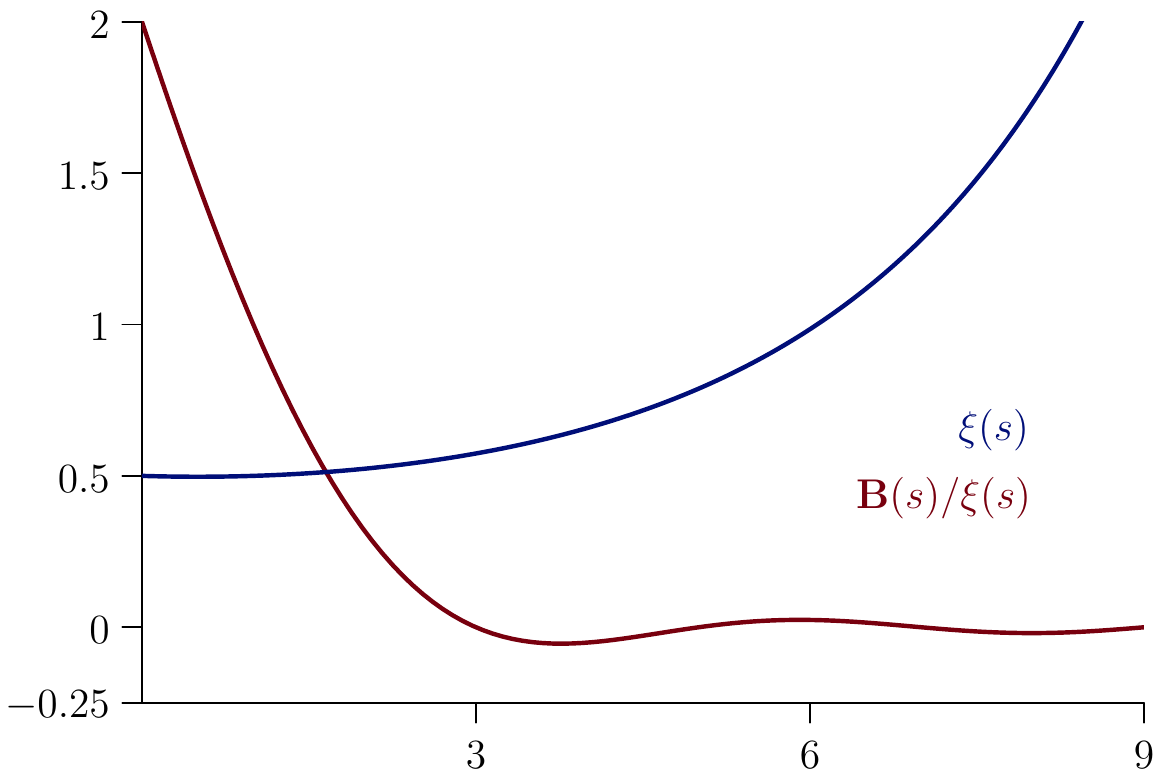}%
    \caption[The Hadamard decomposition]{\textit{The Hadamard decomposition of the Bernoulli function.}}
    \label{fig:fig8hadadecomp}
\end{figure}

\sect{The Hadamard decomposition}{0}

\noindent We denote Hadamard's infinite product over the zeros  of \(\zeta(s) \) by
\begin{equation}
    \operatorname{H}_{\zeta}(s) \, = \,
    \frac12\, \prod_{\Im {\rho} > 0 } \left(1 - \frac{s}{\rho} \right)\left(1 - \frac{s}{1- \rho} \right).
    \label{h99}
\end{equation}
The product runs over the zeros with \( \Im( \rho ) > 0 \).
The absolute convergence of the product is guaranteed as the terms
are taken in pairs \( ( \rho, 1- \rho ) \).
Hadamard's infinite product expansion of \( \zeta(s) \) is
\begin{equation}
    \zeta(s) \, = \, \frac{\pi^{ s/2 } } {\left(s/2\right)! }\, \frac{\operatorname{H}_{\zeta}(s)}{s-1} \, ,
    \quad s \notin \{1, -2, -4,  \dots \}.
    \label{q55}
\end{equation}

Since the zeros of \(\zeta(s) \) and \( \operatorname{B}(s) \) are
identical in the critical strip by (\ref{q52}), this representation
carries directly over to the Bernoulli case. Writing \( \sigma = (1-s)/2 \) we get
(see figure \underline{\ref{fig:fig8hadadecomp}})
\begin{equation}{
        \operatorname{B}(s)\, = \,  \frac{\pi^{\sigma}}{\sigma !} \, \operatorname{H}_{\zeta}(s) \, ,
        \quad s \notin \{-2,-4,-6, \dots \}.
    } \label{q56}
\end{equation}
This is the \textit{Hadamard decomposition of the Bernoulli function.}\index{Hadamard product}

\smallskip

The zeros of \( \operatorname{B}(s) \) with \( \Im( \rho ) = 0 \) are
at \( 3, 5, 7, \ldots \) (making the Bernoulli numbers vanish at these indices),
due to the factorial term in the denominator. This representation separates
the nontrivial zeros on the critical line from the trivial zeros on the real axes.

\smallskip

Here we see one more reason why \(\operatorname{B}_1 = \frac12\).
The oscillating factor has the value \(\frac{\pi^0}{0!} = 1, \)
and the Hadamard factor has the value
\(\operatorname{H}_{\zeta}(1)= - \zeta(0) \cdot 1\).
The Bernoulli value follows from \(\zeta(0) = - \frac12 \).

\smallskip

If we compare the identities \((\ref{q554}) \) and \((\ref{q56}) \), we get as
a corollary \( { \xi =  \operatorname{H}_{\zeta}}. \)
This is precisely the proposition that Hadamard proves in his 1893 paper \cite{Hadamard}.
Thus applying  (\ref{int1}) to (\ref{q55}) leads
to Jensen's formula\index{Jensen!formula} for the Riemann \(\xi\) function,
\begin{equation}
   \xi(s) \, = \, \frac{2 (s/2)!}{\pi^{s/2-1}}
    \int_{- \infty}^\infty\frac{( \tfrac{1}{2} +ix)^{1-s}}{({\Ee}^{\pi x} + {\Ee}^ {- \pi x})^2}\,dx.
    \label{33}
\end{equation}

\sect{The generalized Euler function}{10}

The \textit{generalized Euler function}\index{Euler!generalized function} is defined in terms of
the generalized Genocchi function (\ref{GG}) as
\begin{equation}
    \operatorname{E}(s, v) \, = \, - \frac{\operatorname{G}(s+1, v)}{s+1} . \label{EG}
\end{equation}
Of particular interest are the cases \( v=1\) and \(v=1/2\)
which we will consider now.

\sect{The Euler tangent function}{20}

The  \textit{Euler tangent function}\index{Euler!tangent function} is defined as
\begin{equation}
    \operatorname{E}_{\tau}(s) \, = \, 2^s \operatorname{E}(s, 1),
    \label{eulsimpl}
\end{equation}
where the limiting value \(\log(2)\) closes the definition gap at \(s=-1.\)

An alternative representation is \(\operatorname{E}_{\tau}(s)  = -2 \operatorname{\Re}( {\operatorname{Li}_{-s}(i)} ) \),
where \(\operatorname{Li}_{s}(v) \) denotes the polylogarithm.

Of special interest is  formula (\ref{eulsimpl}) for integer values \(n\),
\begin{equation}
    \operatorname{E}_{\tau}(n) \, = \,   2^n \operatorname{E}(n, 1)  \quad (n \ge 0).
    \label{eultaug}
\end{equation}
These numbers are listed in \OEIS \seqnum{A155585}; the first few values are
\[  \operatorname{E}_{\tau}(n) \, = \, 1,\,1,\,0,\,-2,\,0,\,16,\,0,\,-272,\,  0,\, 7936,\, 0,\, -353792,\,  \ldots \]
Tracing back the definitions and using (\ref{cbrep}),
the Bernoulli function can express the Euler tangent function as
\begin{equation}
    \operatorname{E}_{\tau}(s) \, = \,
    \left(4^{s+1} - 2^{s+1} \right)  \frac{\operatorname{B}(s+ 1)}{s+1}.
    \label{eulbern}
\end{equation}

The Euler tangent numbers are of particular importance because they relate 
the numbers  \(\operatorname{E}_{\tau}(n)\) to another type of Euler numbers, 
the \textit{Eulerian numbers}.\index{Eulerian!number}

\enlargethispage{12pt}
\newpage

\noindent Let \(\genfrac<>{0pt}{}{n}{k}\) denote the Eulerian numbers \OEIS \seqnum{A173018}, then%
\begin{equation}
    \operatorname{E}_{\tau}(n) \, = \,  \sum_{k=0}^n (-1)^k \genfrac<>{0pt}{}{n}{k} .
    \label{eulerian}
\end{equation}
The right side of (\ref{eulerian}) is the value of the \textit{Eulerian polynomial} \cite{LuschnyEL} \index{Eulerian!polynomial}
\(\mathrm{A}_n(x) = \sum_{k=0}^n \genfrac<>{0pt}{}{n}{k} x^k \) at \(x=-1\).
With this we get the identity
\begin{equation}
  2^n \operatorname{E}_n(1) \,  =\,  \operatorname{E}_{\tau}(n)  \, =  \,  \mathrm{A}_n(-1).  
  \label{eulereulerin}
\end{equation}
A further representation results if the factor \(2^n\) is taken into account in the 
Stirling-Fubini representation of the Bernoulli numbers. Indeed, 
\( \operatorname{E}_{\tau}(n) = \operatorname{P}_n(1) \) where 
\begin{equation}
   \operatorname{P}_n(x) = \sum_{k=0}^n (-2)^{n - k}  \stirling{n}{k} k! \, x^k .
     \label{stirfubinieul}
\end{equation}

\sect{The Euler secant function}{10}

\begin{figure}
    \centering
    \includegraphics[height=240pt,width=320pt]{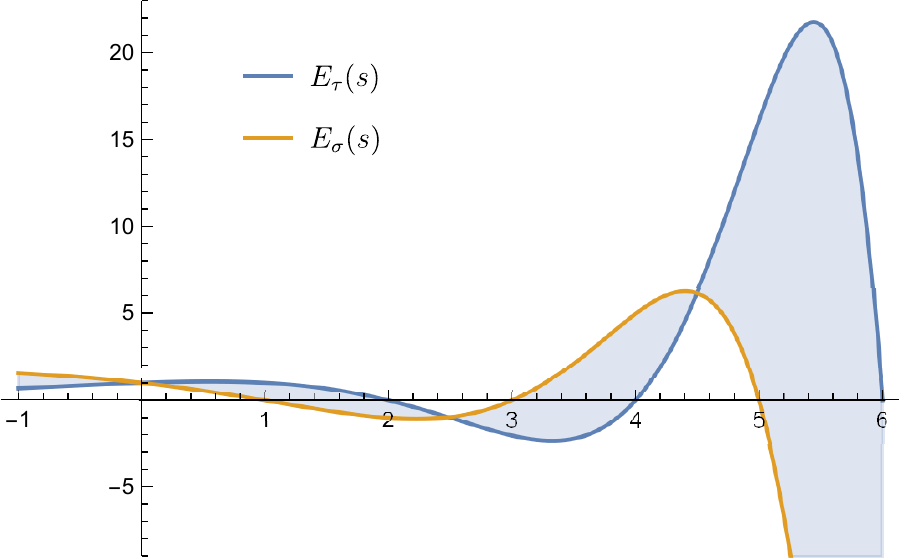}%
    \caption[Euler sec vs. tan]{\textit{Euler tangent versus Euler secant function.}} 
    \label{fig:fig23eulersectan}
\end{figure}

\noindent The  \textit{Euler secant function}\index{Euler!secant function} is defined
\begin{equation}
    \operatorname{E}_{\sigma}(s) \, = \, 2^{s} \operatorname{E}\left(s,  1 / 2 \right),
    \label{EulerCentral}
\end{equation}
with the limiting value \( \operatorname{E}_{\sigma}(-1)  \, = \, \pi/2. \)  

\enlargethispage{12pt}
\newpage

We also write \( \operatorname{E}(s) = \operatorname{E}_{\sigma}(s)\) and call 
this function the \textit{Euler function}\index{Euler!function}, because it interpolates the numbers
\begin{equation}\operatorname{E}_n \ = \  \operatorname{E}_{\sigma}(n)
    \ = \ 2^{n} \operatorname{E}\left(n,  1 / 2 \right)  \quad (n \ge 0),%
    \label{EulerNum}
\end{equation}
which traditionally are called the \textit{Euler numbers}.\index{Euler!numbers} 

\smallskip
The \OEIS identifier is \seqnum{A122045} and the sequence starts
 \[  \operatorname{E}_n \, =  1,\, 0,\, -1,\, 0,\, 5,\, 0,\, -61,\, 0,\, 1385,\, 0,\, -50521,\, 0,\, 2702765,\,  \ldots \]

The generalized Bernoulli function can represent the Euler secant function as
\begin{equation}
    \operatorname{E}_{\sigma}(s) \, = \,
   \frac{2}{s+1}  \,  4^{s} \left( \operatorname{B} \left(s+1, \frac34 \right) - \operatorname{B} \left(s+1, \frac14 \right) \right)  .
    \label{eulnumber}
\end{equation}
A more compact form is
\(\operatorname{E}_{\sigma}(s)  = 2 \operatorname{\Im}( \operatorname{Li}_{-s}(i) )  \).

\smallskip
The \textit{Jensen formula} for the even-indexed classical Euler numbers follows
from (\ref{eulnumber}) and the Jensen formula (\ref{int5}) for the generalized
Bernoulli function.
\begin{equation}
    \operatorname{E}_{2n} \, = \,  \frac{2 \pi}{2n+1}
    \int_{0}^\infty   \frac{ (4zi+1)^{2n+1} - (4zi-1)^{2n+1} }
    {({\Ee}^{- \pi z} + {\Ee}^ {\pi z})^2}\, dz.
    \label{jenseul}
\end{equation}

\sect{The family of Bernoulli and Euler numbers}{20}

\noindent The family tree of the Euler numbers is subdivided into three branches:
the \textit{Euler secant numbers}, the \textit{Euler tangent numbers},
and the \textit{Euler zeta numbers}. The Euler zeta numbers are rational numbers,
whereas the Euler secant and Euler tangent numbers are integers.

\smallskip
The traditional way of naming reserves the name `Euler numbers' for the
Euler secant numbers, while the way preferred by combinatorialists
(see for instance Stanley \cite{AltPerm}) is to call Euler numbers
what we call the Andr\'e numbers.

\smallskip
We introduce the name `Andr\'e numbers' in honor of D\'esir\'e Andr\'e, who
studied their combinatorial interpretation as \(2\)-alternating permutations
in 1879 and 1881 \cite{Andre79, Andre81}.
We believe that this is a fair sharing of the mathematical  name space and eliminates
the ambiguity that otherwise exists.

\newpage
\begin{figure}[ht]
    \centering
    \includegraphics[width=110mm]{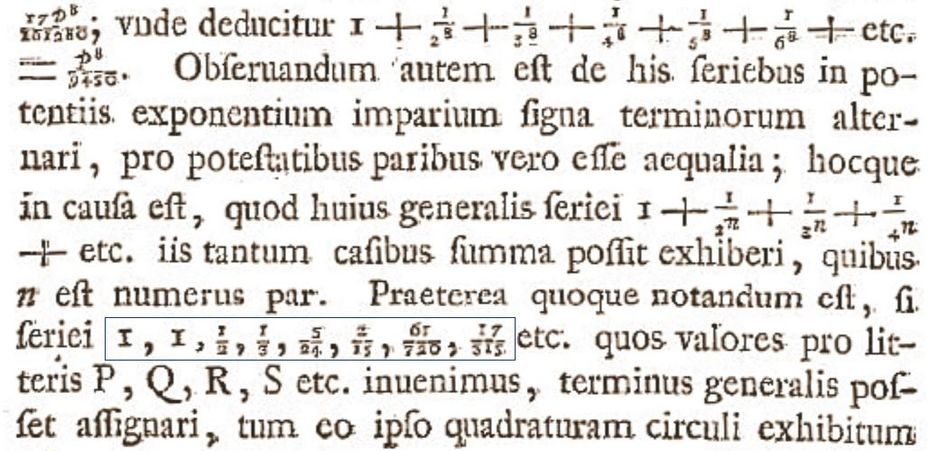}%
    \caption[L. Euler: De summis serierum]{\textit{L. Euler, De summis serierum reciprocarum, 1735.}}%
\end{figure}

 \bigskip
{\centering
    \renewcommand{\arraystretch}{1.3}
    \begin{tabular}{c|cccccccc}
        \rule[-1ex]{0pt}{2.5ex} \(n\)                                                 & 1      & 2      & 3         & 4     & 5       & 6          & 7              & 8      \\     \hline
        \rule[-1ex]{0pt}{2.5ex} \(\operatorname{E}_{\sigma}\) & $0$ & $-1$ & $\ 0$  & $5$ & $0$  & $-61$   & $\ \ 0$   & $1385$ \\
        \rule[-1ex]{0pt}{2.5ex} \(\operatorname{E}_{\tau}\)    & $1$  & $0$  & $-2$    & $0$ & $16$ & $0$      & $-272$    & $0$    \\
        \rule[-1ex]{0pt}{2.5ex} \(\mathcal{A}  \)  & $1$  & $1$  & $\ \ 2$ & $5$ & $16$ & $\ \ 61$ & $\ \ 272$ & $1385$ \\
    \end{tabular}
    \captionsetup{hypcap=false}   \captionof{table}{\textit{Euler and Andr\'e numbers}}}

\bigskip
Let us try to replicate the above extension procedure for the Bernoulli numbers.
The next table shows the outcome of our choice.

\bigskip
{ \centering
    \renewcommand{\arraystretch}{1.3}
    \begin{tabular}{c|cccccccc}
        \rule[-1ex]{0pt}{2.5ex} \(n\)                         & 1         & 2             & 3                   & 4                  & 5                & 6              & 7                       & 8    \\   \hline
        \rule[-1ex]{0pt}{2.5ex} \(\operatorname{B}_{\tau}\)   & $0$       & $\frac{1}{6}$ & $0$                 & $- \frac{1}{30}$   & $0$              & $\frac{1}{42}$ & $0$                     & $- \frac{1}{30} $  \\
        \rule[-1ex]{0pt}{2.5ex} \(\operatorname{B}_{\sigma}\) & $\frac12$ & $\ 0$         & $- \frac{3}{56}$    & $\ \ 0$            & $\frac{25}{992}$ & $0$            & $- \frac{427}{16256}$   & $\ \  0$           \\
        \rule[-1ex]{0pt}{2.5ex} \(\mathcal{S}\)          & $\frac12$ & $\frac16$     & $ \ \ \frac{3}{56}$ & $\ \ \frac{1}{30}$ & $\frac{25}{992}$ & $\frac{1}{42}$ & $\ \ \frac{427}{16256}$ & $\ \ \frac{1}{30}$ \\
    \end{tabular}
    \captionsetup{hypcap=false}   \captionof{table}{\textit{Bernoulli and Seki numbers}}
    \label{berntable2}
}

\bigskip
\noindent We introduce the name `Seki numbers' \index{Seki!numbers} in honor of Takakazu Seki, who
discovered Bernoulli numbers before Jacob Bernoulli (see \cite{Seki}) to denote
the extended version of the Bernoulli numbers in their unsigned form, 
which is the third row in the table above. 

The relationship between the seven sequences is shown in figure \underline{\ref{fig:EulerBernoulli}}.
It reveals  that the Bernoulli numbers and the Euler numbers have a common
backbone: the \textit{Euler zeta numbers}\index{Euler!zeta numbers}.
Euler introduced these rational numbers in 1735 in
\textit{De summis serierum reciprocarum} \cite{EulerE41Eng}.
We will denote the numbers with \( \mathcal{Z}_n \). They begin for  \(n \ge 0 \) %
\begin{equation}
    1, \ 1,\  \frac12,\  \frac{1}{3},\ \frac{5}{24},\ \frac{2}{15},\  \frac{61}{720},\  \frac{17}{315}, \ \frac{277}{8064},\ \frac{62}{2835}, \ldots .
    \label{EulerZetaNum}
\end{equation}

\newpage

\begin{sidewaysfigure}
    \centering
    \includegraphics[width=210mm]{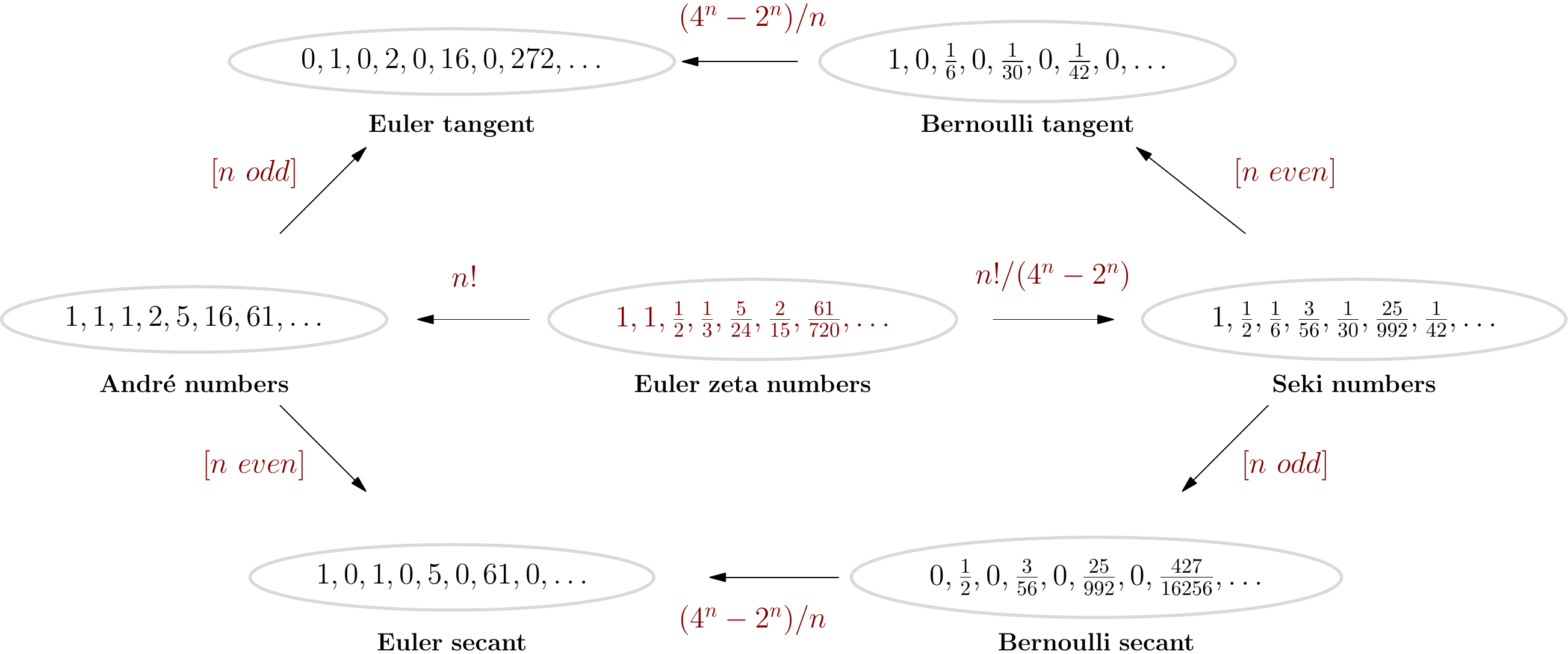}
    \caption[The Euler--Bernoulli number family]
    {\textit{ \phantom{.} \\The Euler--Bernoulli family of numbers (unsigned version)}}
    \label{fig:EulerBernoulli}
\end{sidewaysfigure}

\FloatBarrier \newpage

\sect{The Bernoulli secant numbers}{0}

\noindent The \textit{Bernoulli secant function}\index{Bernoulli!secant function}
\( \operatorname{B}_{\sigma} \) is defined by the equation 
\begin{equation}
    \operatorname{B}_{\sigma}(s) \, = \,  \frac{2^{s - 1}}{2^s - 1}
    \left (\operatorname{B}\left(s, \frac34\right) - \operatorname{B}\left(s, \frac14\right) \right) .           
    \label{bernsec}    
\end{equation} 
\textit{Bernoulli secant numbers}\index{Bernoulli!secant numbers} are for integer \(n \ge 1\) the values of 
the Bernoulli secant function and by convention \(\operatorname{B}^{\sigma}_0 =0\).
\begin{equation}
    \operatorname{B}^{\sigma}_n  \, = \,   0,\  \frac12,\  0,\   -\frac{3}{56},\  0,\  \frac{25}{992},\  0,\  -\frac{427}{16256},\  0 , \ \ldots
\end{equation} 
Thus, if one calls the classical Bernoulli numbers the 
\textit{Bernoulli tangent numbers}\index{Bernoulli!tangent numbers} \(\operatorname{B}^{\tau}_n\), 
one gets a way of speaking that corresponds to the classical terminology associated
with the Euler numbers. See figure \underline{\ref{fig:EulerBernoulli}}, and \OEIS \seqnum{A160143}, \seqnum{A193476}.

\smallskip
The Bernoulli secant numbers can be  represented by the Euler secant numbers since 
\begin{equation}
    \operatorname{B}_{\sigma}(n)\,=\, (-1)^{n - 1} \frac{n}{4^n - 2^n} \operatorname{E}_{\sigma}(n-1) \quad (n \ge 1) .
    \label{berneulsec}
\end{equation}
From (\ref{eulnumber}) and (\ref{bernsec}) we see that this is a special case of
\begin{equation}
    \operatorname{B}_{\sigma}(s) \, = \,   \frac{2 s}{4^s - 2^s} \operatorname{\Im}( \operatorname{Li}_{1 - s}(i) ) .
\end{equation}

\begin{figure}
    \centering
    \includegraphics[height=240pt,width=320pt]{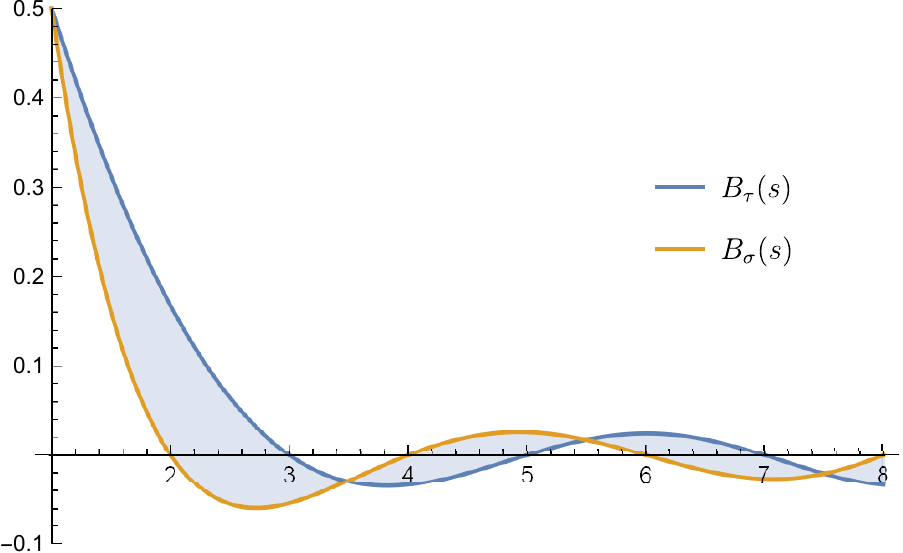}%
    \caption[Bernoulli sec vs. tan]{\textit{Bernoulli tangent versus Bernoulli secant function.}} 
    \label{fig:fig24bernoullisectan}
\end{figure}

\sect{The extended Bernoulli function}{20}

\noindent  The \textit{extended zeta function} \index{Riemann!extended function} is defined for \(s \ne 1\) as
\begin{equation} {
    \widetilde{\zeta}(s)  \, = \,   \zeta \left({s} \right) - \frac{ \zeta \left({s},{\frac{3}{4}} \right)
    - \zeta \left({s},{\frac{1}{4}} \right) } {{2}^{{s}}-2} .
    \label{lzeta}
    } \end{equation}
The \textit{extended Bernoulli function}\index{Bernoulli!extended function}
is defined for \(s \ne 0\) as
\begin{equation} {
        \mathcal{B}(s) \, = \, -s\, \widetilde{\zeta}(1-s).
    } \label{genbern}
\end{equation}
Setting \(\mathcal{B}(0) = 1 + {\pi}/{\log(4)}\) closes the definition gap.

\noindent An equivalent  definition using the generalized Bernoulli function is%
\begin{equation}
    \mathcal{B}(s) \, = \, \operatorname{B}(s)  +  \frac{2^{s - 1}}{2^s - 1}
    \left (\operatorname{B}\left(s, \frac34\right) - \operatorname{B}\left(s, \frac14\right) \right) .
    \label{BernExt}
\end{equation}
With the terms  introduced above, this says that the extended Bernoulli function is the
the sum of the tangent Bernoulli and the secant Bernoulli function.
\begin{equation}
    \mathcal{B}(s) \, = \operatorname{B}_{\tau}(s) + \operatorname{B}_{\sigma}(s). 
\end{equation}

The \textit{extended Bernoulli numbers}\index{Bernoulli!extended numbers} are the values of the
extended Bernoulli function at the positive integers. 
\begin{equation}
    1,\   \frac16,\   -\frac{3}{56},\  -\frac{1}{30},\  \frac{25}{992},\  \frac{1}{42},\  -\frac{427}{16256},\  -\frac{1}{30} , \ \ldots
\end{equation}
  
\sect{The extended Euler function}{14}

\noindent The \textit{extended Euler function}\index{Euler!extended function}
is a generalization of the identity (\ref{eulbern}),
\begin{align}
    \mathcal{E}(s) \, & = \, \left( 4^{s+1} - 2^{s+1} \right)
    \frac{\mathcal{B}(s+1)}{ (s+1)},        \label{exteulbern}                           \\
                      & = \,  (2^{s+1} - 4^{s+1})\, \widetilde{\zeta}(-s).
    \label{exteulzeta}
\end{align}
The limiting value \( \pi/2 + \log(2) \) closes the definition gap at \(s=-1\).%

\newpage

The extended Euler function is the sum of  the Euler secant and the Euler tangent function.
\begin{equation}
    \mathcal{E}(s)    \, = \, \operatorname{E}_{\sigma}(s) + \operatorname{E}_{\tau}(s)
    \, = \, 2^s (\operatorname{E}\left(s,  1/2 \right) + \operatorname{E}\left(s,  1 \right) ).
    \label{anum}
\end{equation}

\noindent For integer \(n \ge 0\)  we write \( \mathcal{E}_n  =    \mathcal{E}(n) \). The sequence starts:
\begin{equation}
    \mathcal{E}_n  =   2,\, 1,\, -1, \, -2,\,  5,\, 16,\,-61,\, -272,\, 1385,\, \dots.
    \label {extenum}
\end{equation}
These are the  \textit{extended Euler numbers}  \(\mathcal{E}_n\), \OEIS \seqnum{A163982} negated.

\smallskip
Since \(\operatorname{B}(s, 1/2) = \operatorname{B}(s) (2^{1-s} - 1)\),
the difference of the right-hand sides of
(\ref{eulnumber}) and  (\ref{eulbern})  reduces to
\begin{equation}
    \mathcal{E}(s) \, =\, \frac{ (4^{s+1} - 2^{s+2}+2) \operatorname{B}(s+ 1)
        - 4^{s+1}  \operatorname{B} (s+1, 1/4) }  {s+1}.
    \label{andnum}
\end{equation}
The \textit{Jensen representation}  follows from (\ref{andnum}).
\begin{equation}
    \mathcal{E}(s) \, = \,  \frac{2 \pi}{s+1}
    \int_{\infty}^\infty
    \frac{ (4 z i + 1)^{s + 1} + (1 - 2^{-2s - 1})(2 + 4 z i)^{s + 1} }
    {({\Ee}^{- \pi z} + {\Ee}^ {\pi z})^2}\, dz.
    \label{jensand}
\end{equation}

\sect{The Andr\'e function}{10}

\noindent Many of the traditional integer sequences considered here are 
signed, like the Bernoulli numbers and the Euler numbers. 
However, under the influence of combinatorics, more and more the 
unsigned versions of these numbers have come into focus. 
The paradigmatic example are the absolute Euler numbers, 
which we call Andr\'e numbers.

\smallskip
The unsigned versions of the Euler secant and Euler tangent functions are defined as
\index{Euler function!abs secant}\index{Euler function!abs tangent} 
\begin{gather} 
        \operatorname{|E|}_{\sigma}(s) \, = \, \cos(\pi s / 2) \operatorname{E}_{\sigma}(s); \\
        \operatorname{|E|}_{\tau}(s) \, = \, \sin(\pi s / 2) \operatorname{E}_{\tau}(s).
        \label{UnsignedEulerFunctions}
\end{gather}
The \textit{Andr\'e function}\index{Andr\'e! function} is defined as the sum of these two 
unsigned Euler functions,
\begin{equation} {
        \mathcal{A} (s) \, = \, \operatorname{|E|}_{\sigma}(s) + \operatorname{|E|}_{\tau}(s).
        \label{UnsignedAndreFun}
}    \end{equation}
Considering the long chain of definitions on which (\ref{UnsignedAndreFun}) is based, 
it is astonishing how easily it can be represented by a single function.%
\begin{equation} {
        \mathcal{A} (s) \, = \, (-i)^{s+1} \operatorname{Li}_{-s}(i) + i^{s+1} \operatorname{Li}_{-s}(-i)  .
        \label{UnsignedAndreFunctionPL}
}    \end{equation}
Here \(i\) is the imaginary unit, \(\operatorname{Li}_{s}(v)\) is the polylogarithm,
and the principal branch of the logarithm is used for the powers.

\begin{figure}
    \centering
    \includegraphics[height=240pt,width=320pt]{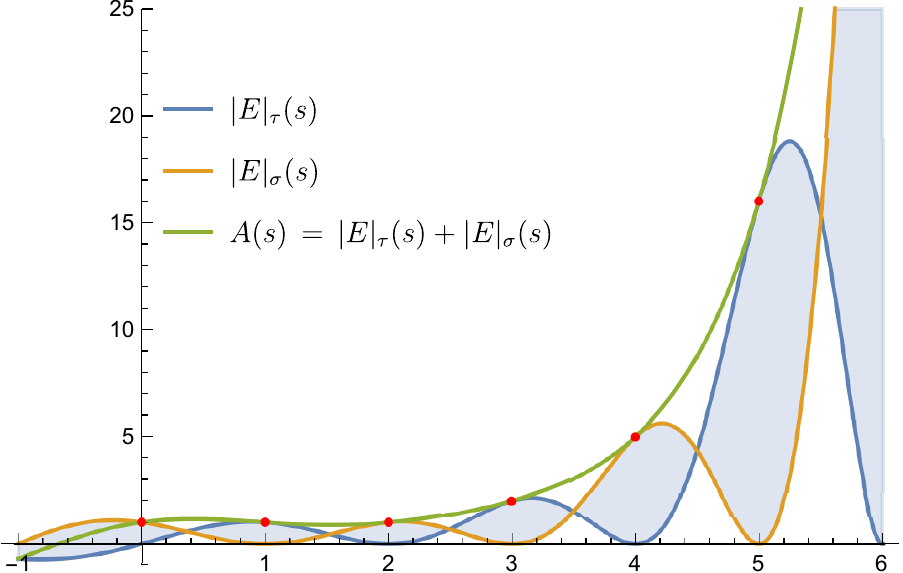}%
    \medskip
    \caption[Andre function]{\textit{The Andr\'e function} \(\mathcal{A}(s).\)}
    \index{Andre!function}
    \label{fig:fig37andrefunction}
\end{figure}    
\FloatBarrier

The \textit{Andr\'e numbers}\index{Andr\'e!numbers}
\(\mathcal{A}_n = \mathcal{A} (n)\) are listed
for integer \(n \ge 0\) in \OEIS \seqnum{A000111}.
The sequence starts:
\begin{equation}
    \mathcal{A}_n  = 1,\, 1,\, 1, \, 2,\,  5,\, 16,\, 61,\, 272,\, 1385,\, 7936, \, \dots.
\end{equation}
For integer \(n \ge 1\) equation (\ref{UnsignedAndreFunctionPL}) simplifies to
\begin{equation}
\mathcal{A}_n  = 2\, (-i)^{n+1}  \operatorname{Li}_{- n}(i).%
\end{equation}
The \textit{Euler zeta numbers} are, by definition, 
\begin{equation} {
        \mathcal{Z}_n \, = \frac{\mathcal{A}_n}{n !}  \quad (n \ge 0).
        \label{eulzetanum}
} \end{equation}
The \textit{signed Andr\'e function}\index{Andr\'e!signed function} \( \mathcal{A}^{\ast}(s) \) 
interpolates the \textit{signed Andr\'e numbers} \index{Andr\'e!signed numbers}  
\( \mathcal{A}^{\ast}_n = (-1)^{n} \mathcal{A}_n, \) and is defined as
\begin{equation} 
    \mathcal{A}^{\ast}(s) \, =\, i e^{-\frac{1}{2} i \pi  s} \left(\text{Li}_{-s}(-i)-e^{i \pi  s} \text{Li}_{-s}(i)\right). 
    \label{SignedAndreFun}
 \end{equation} 
The signed Andr\'e numbers  are \OEIS \seqnum{A346838},
and differ from (\ref{extenum}) in the first term and by the sign pattern.
\begin{equation}
    \mathcal{A}^{\ast}_n  =   1,\, -1,\, 1, \, -2,\,  5,\, -16,\, 61,\, -272,\, 1385,\, \dots.
    \label{sigandnum}
\end{equation}

\sect{The Seki function}{0}

\begin{figure}    
    \centering
    \includegraphics[height=240pt,width=320pt]{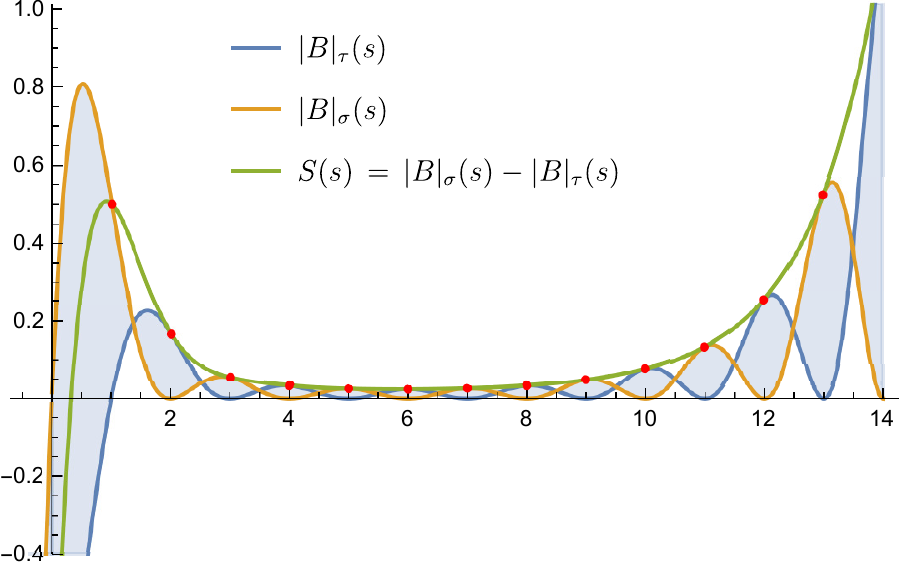}%
    \medskip
    \caption[Seki function]{\textit{The Seki function} \(\mathcal{S}(s).\)}
    \index{Seki!function}
    \label{fig:fig38sekifunction}
\end{figure}

\noindent The unsigned versions of the Bernoulli secant and Bernoulli tangent functions are defined as
\index{Bernoulli function!abs secant}\index{Bernoulli function!abs tangent} 
\begin{gather} 
    \operatorname{|B|}_{\sigma}(s) \, = \, \sin(\pi s / 2) \operatorname{B}_{\sigma}(s); \\
    \operatorname{|B|}_{\tau}(s) \, = \, \cos(\pi s / 2) \operatorname{B}_{\tau}(s).
    \label{UnsignedSekiFunctions}
\end{gather}
The \textit{Seki function}\index{Seki!function} is defined as the difference of these two 
unsigned Bernoulli functions,
\begin{equation} {
        \mathcal{S} (s) \, = \, \operatorname{|B|}_{\sigma}(s) - \operatorname{|B|}_{\tau}(s).
        \label{SekiFun}
}    \end{equation}
One can also base the Seki function  on the Andr\'e function \( \mathcal{A}\).
For \( s=0 \) we set  the limiting value, \( \mathcal{S} (0) = -1 \), and otherwise 
\begin{equation} {
        \mathcal{S} (s) \, = \, \frac{s\, \mathcal{A} (s-1)}{4^s - 2^s}  \quad (s \ne 0).
        \label{SekiAndreFun}
}    \end{equation}
In terms of the polylogarithm this is
\begin{equation} {
        \mathcal{S} (s) \, = \, \frac{s}{4^s - 2^s}\left( (-i)^{s} \operatorname{Li}_{1-s}(i) 
                                 + i^{s} \operatorname{Li}_{1-s}(-i) \right) \quad (s \ne 0).
        \label{SekiFunPL}
}    \end{equation}
The Seki function \({\mathcal{S}}\) interpolates the absolute values of the extended
Bernoulli numbers for \(n \ge 1\). The \textit{Seki numbers}\index{Seki!numbers}
are \(\mathcal{S}_n = \mathcal{S}(n) \) if \(n \ge 1\), and \(\mathcal{S}_0 = 1\) by convention.

\smallskip
The \textit{signed Seki function}\index{Seki!signed function} \( \mathcal{S}^{\ast}(s) \) 
interpolates the \textit{signed Seki numbers} \index{Seki!signed numbers}  
\( \mathcal{S}^{\ast}_n = (-1)^{n} \mathcal{S}_n, \) and is defined as \( \mathcal{S}^{\ast}(0) = 1 \), and otherwise%
\begin{equation} 
    \mathcal{S}^{\ast}(s) \, =\, \frac{s}{2^s - 4^s}\,  e^{i \pi s / 2} \left( e^{i \pi s} \text{Li}_{1 - s}(-i) + \text{Li}_{1 - s}(i) \right).
    \label{SignedSekiFun}
\end{equation} 
The signed Seki numbers  differ from (\ref{extenum}) in the first two terms and by the sign pattern.
Apart from the signs, this is the sequence \OEIS \seqnum{A193472}/\seqnum{A193473}.

\sect{The Swiss-knife polynomials}{20}

\noindent The Euler equivalent to the central Bernoulli polynomials is the sequence of \textit{Swiss-knife polynomials}
\index{Swiss-knife!polynomials}\index{Polynomials!Swiss-knife}(see figure \underline{\ref{fig:fig20swissknife}}),
defined as the Appell sequence associated with the Euler numbers \(\operatorname{E}_n\).
\begin{equation}
    \operatorname{K}_n(x) \, = \,  \sum_{k=0}^n \binom{n}{k} \operatorname{E}_k x^{n-k}.
    \label{swiss}
\end{equation}
\begin{table}[b]
    \centering \setlength{\extrarowheight}{3pt}
    {  \begin{tabular}{c}%
            \(1\)                                                                                           \\
            \(x\)                                                                                            \\
            \({x}^{2}\,- \,1\)                                                                         \\
            \({x}^{3}\,- \, 3\,x\)                                                                    \\
            \({x}^{4}\,- \, 6\,x^2\,+\,5\)                                                     \\
            \(x^5\,- \,10\,x^3\,+\,25\,x\)                                                   \\
            \(x^6\,- \,15\,x^4\,+\,75\,x^2\,- \,61\)                                   \\
            \(x^7\,- \,21\,x^5\,+\,175\,x^3\,- \,427\,x\)                           \\
            \(x^8\,- \,28\,x^6\,+\,350\,x^4\,- \,1708\,x^2\,+\,1385\)    \\
            \(x^9\,- \,36\,x^7\,+\,630\,x^5\,- \,5124\,x^3\,+\,12465\,x\) \\
    \end{tabular} }
    \caption{\textit{The Swiss-knife polynomials} \(\Kappa_n(x)\) for \(n \ge 0.\)}
    \label{skppoly}
\end{table}
They were introduced  in \OEIS \seqnum{A153641} and \seqnum{A081658}.
The author discussed them in \cite{LuschnySKP} and dubbed them \textit{Swiss-knife polynomials}
\(\operatorname{K}_n(x)\) because they allow calculating the Euler--Bernoulli family of
numbers efficiently. The coefficients of the polynomials are integers, in contrast to
the coefficients of the Euler and Bernoulli polynomials.
The parity of the monic  \(\operatorname{K}_n(x)\) equals the parity of \(n\).

The \textit{Worpitzky representation} of the Swiss-knife polynomials is
the generalized Worpitzky transform (\ref{q46})
of the sequence \(\nu(n)\), where \(\nu(n) = \cos((n-1)\pi/4)/2^{(n-1)/2} \).

\enlargethispage{12pt}
\newpage

   \begin{figure}
    \centering
    \includegraphics[height=170pt,width=310pt]{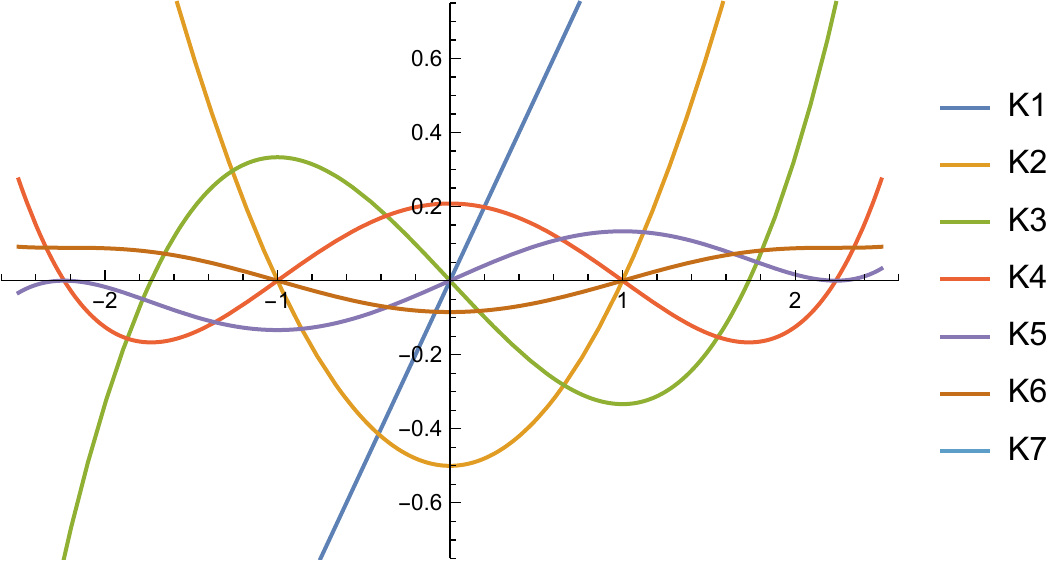}%
    \caption[The Swiss knife polynomials]{\textit{The Swiss knife polynomials} \(\Kappa_n(x)/n!\)}
     \label{fig:fig20swissknife}
\end{figure}

Equivalently, let \(\alpha\) be the repeating sequence with period \(p\), where
\( p = (1, 1, 1, 0, -1, -1, -1, 0)\), then
\begin{equation}{
        \operatorname{K}_n(x) = \sum_{k=0}^{n} \frac{\alpha(k)}{2^{ \lfloor k/2
                    \rfloor}} \sum_{v=0}^{k}(-1)^{v} \binom{k}{v}(x+v+1)^n.
        \label{skpworp}
    } \end{equation}

Chen \cite[theorem 3]{Chen} proves (\ref{skpworp}) for the Euler numbers \( \operatorname{K}_n(0)\) and
\(\operatorname{K}_n(1)\) using the Akiyama--Tanigawa algorithm.

The Swiss-knife polynomials can be computed efficiently with the following \textit{recurrence}:

\smallskip
\noindent\rule{\textwidth}{0.4pt} \newline 

\indent Set \( \operatorname{K}_{0}(x) = 1 \)  for all \(x\).

Now assume \( \operatorname{K}_{n-1}(x) \) already computed and take the coefficients%
\[ 
    c_k^{(n-1)}  = [x^{n-2k}] \operatorname{K}_{n-1}(x) \text{ for } k \in \{0,1,\ldots, \lfloor (n-1)/2 \rfloor \}. 
\]

Next compute
\begin{equation}
    c_k^{(n)} =   c_k^{(n-1)} n / (n-2k)  \ \ \text{for } k \in \{0,1,\ldots, \lfloor (n-1)/2 \rfloor \} .
\end{equation}
\indent If \(n\) is even, set additionally  \( c_{n/2}^{(n)} = -\sum_{k=0}^{\lfloor (n-1)/2 \rfloor} c_k^{(n)} \). Then%
\begin{equation}
    \operatorname{K}_{n}(x) = \sum_{k=0}^{ \lfloor n/2 \rfloor}  c_k^{(n)} x^{n-2k}.
\end{equation}

\medskip
\noindent\rule{\textwidth}{0.4pt}

\medskip
The algorithm is based on the fact that the Swiss-knife polynomials form an Appell sequence.

\enlargethispage{12pt}
\newpage

\begin{figure}
    \centering  \index{Seki!approximation}
    \includegraphics[height=240pt,width=320pt]{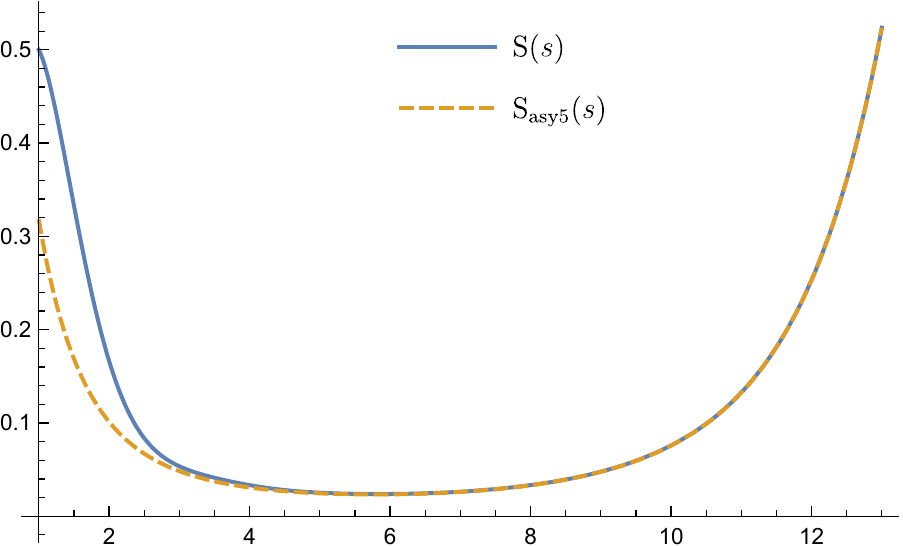}%
    \medskip
    \caption[Euler Zeta numbers]{\textit{Approximation to the Seki function.}}
    \caption[]{\( \mathcal{S}(s) \, \simeq \, 2^{\frac{3}{2}-s} \pi ^{\frac{1}{2}-s} s^{s+\frac{1}{2}} e^{\frac{(2/7) - s^2 + 30 s^4 - 360 s^6} {360 s^5}} \)}
    \label{fig:fig33approxabsbern}
\end{figure}

\sect{Asymptotics for the Bernoulli function}{0}

\noindent An asymptotic expansion of the Bernoulli function follows directly from
\(\operatorname{B}(s) = -  \zeta(s)\, \tau(s)\, s! \) by using
Stirling's formula and the generalized harmonic numbers.

For the remainder term we will use the notation
\begin{equation}
    \operatorname{R_K}(s) \, = \,  \exp\left(\frac12 +  \sum_{n = 1}^{\operatorname{K}}
    \frac{\operatorname{B}(n+1)}{n+1} \frac{s^{-n}}{n} \right).
    \label{asyrem}
\end{equation}

For an \textit{even positive integer} \(n\), the Bernoulli function has an
efficient asymptotic approximation\index{Bernoulli!asymptotics} \cite{Luschny5},
which reads with \(\operatorname{R_5}(s)\)
\begin{equation} {
        |\operatorname{B}(n)|  \,  \sim \, 4\,\pi\, \left( \frac {n}{2 \pi \Ee} \right)^{n+1/2}
        \exp\left(\frac12 + \frac{n^{-1}}{12} - \frac{n^{-3}}{360} + \frac{n^{-5}}{1260}\right) .}
    \label{bernumasy}
\end{equation}
The coefficients originate from the Stirling expansion of  \( \log(\Gamma(s)) \)
(see \seqnum{A046969}). 

The number of exact decimal digits guaranteed
by (\ref{bernumasy}) is \( 3 \log(3\, n) \) if \(n \ge 50\).
The Boost \( \operatorname{C}^{++} \) library \cite{Boost} uses
this approximation  for huge arguments \(n\).

Different asymptotic developments can be based on other expansions of
the Gamma function, for instance, on Binet's formula \cite[p.\ 48, \seqnum{A122252}]{Erdelyi}
generalized by Nemes \cite[4.2]{Nemes13}. More general asymptotic
expansions and error bounds follow from those of the Hurwitz zeta
function established in Nemes \cite{Nemes}. 

\smallskip
But much more is true: The assumption made for (\ref{bernumasy}) that \(n\) is an even positive
integer can be dropped  if one adds the factor \(- \cos \left( \frac{s \pi}{2} \right)\) to the right side.
This gives the asymptotic expansion of the Bernoulli function for
real \(s > 0\), with \(\operatorname{K}\) suitably chosen,%
\begin{equation}
    \operatorname{B} (s)  \, \sim \,   4\, \pi \,
    \left({\frac { s }{{{2 \pi \Ee}} }}\right)^{s+1/2}
    \left( -  \cos \left( \frac{s \pi}{2} \right) \right)
    \operatorname{R}_{\operatorname{K}}(s).
    \label{asygen}
\end{equation}
Moreover, the Seki \(\mathcal{S}(s)\) has the
corresponding expansion without the factor \(- \cos \left( \frac{s \pi}{2} \right)\).

\smallskip
The close connection between the Bernoulli and the Euler numbers 
is also reflected in the fact that the asymptotic development 
of the Euler function differs formally only slightly from (\ref{asygen}).
\begin{equation}
    \operatorname{E} (s)  \, \simeq \,  - 4\, 
    \left({\frac {2  s }{{{ \pi \Ee}} }}\right)^{s+1/2}
    \left( -  \cos \left( \frac{s \pi}{2} \right) \right)
    \operatorname{R}_{\operatorname{K}}(s).
    \label{asyeul}
\end{equation}

From (\ref{asygen}) and (\ref{asyeul})  asymptotic expansions for other functions can be 
easily derived.  As an example we show an asymptotic expansion of the logarithm of 
the Andr\'e function of order \(O(s^{-7})\).
\begin{equation}
    \log \mathcal{A}(s) \  \simeq \ \log(4) +  \left(\frac12 + s \right) \log\left(\frac{2 s}{\pi}\right)
    +\frac{(2/7) - s^2 + 30 s^4 - 360 s^6}{360 s^5}
\end{equation}

\sect{Definition dependencies}{0}

\noindent In the present essay we have entirely dispensed with generating
functions and have only taken the analytical point of view. This resulted in a 
net of hierarchically structured definitions shown in the graph \underline{\ref{defdepend}} below.
The numbers attached to the arrows indicate the corresponding formula.

This shows that all functions of the Euler-Bernoulli family can be represented 
using the central formula (\ref {int5}). This approach is not just theoretically interesting; 
it might also offer computational advantages.

To this end, we note that integrals of the Jensen type  can be evaluated
numerically efficiently to high accuracy.
Based on Johansson and Blagouchine \cite{Johansson}
relevant routines were implemented by Johansson in an arbitrary-precision 
software library \cite{JohanssonArb} with rigorous bounds and used to compute the Stieltjes constants. 

Further expansion of this computational infrastructure to the Hurwitz-Bernoulli function
providing an alternative to the Hurwitz-Riemann function would be highly desirable.

\begin{figure}[t]
    \hspace{-8pt} 
    \includegraphics[height=240pt,width=320pt]{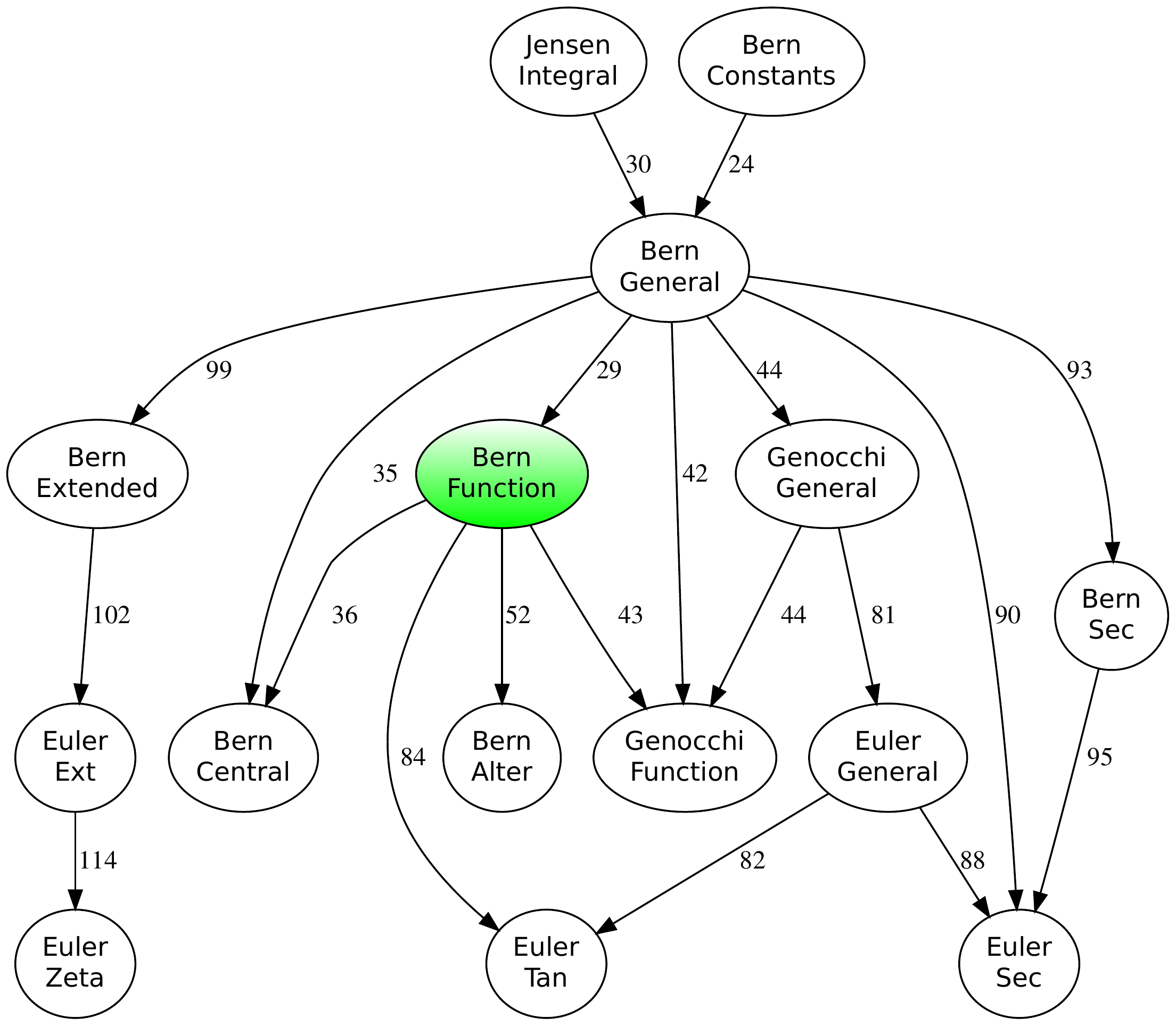}%
    \medskip
    \caption[Dependency graph of definitions]{Dependency graph of definitions.}%
    \label{defdepend}
\end{figure}

\sect{Epilogue: Generating functions}{20}

\noindent The value of \( \operatorname{B}(1)\) deserves special attention. Since it is well known
that \( \sum_{j=0}^\infty {\gamma_j}/ {j!} = 1/2 \), it follows from (\ref{q7}) that
\( \operatorname{B}(1) = 1/2\). Unfortunately, the popular generating
function \( z/(\mathrm{e}^{z}-1) \)\index{generating function!inconsistent} misses this value
and disrupts at this point the connection between the Bernoulli numbers
and the \( \zeta \) function.

\smallskip
For those who do not care about the connection with the zeta function, we add:
Even the most elementary relations between
the Bernoulli numbers and the Bernoulli polynomials break with this choice.
For instance, consider the basic identity (\ref{bernoulligoeshalf}).
It applies to all Bernoulli numbers if
\(\operatorname{B}_n = \operatorname{B}_n(1)\)
but not if \(\operatorname{B}_n = \operatorname{B}_n(0)\) is set.

\smallskip
Many other relations get restricted in their range of validity if the wrong choice is made,
for example, the relation between the Bernoulli numbers and the Eulerian numbers.
Such examples are described in the discussion \cite{Luschny}.

Instead, use the power series \(f(z)\) with the constant term \(1\) such that the
coefficient of \(x^n\) in \((f(x))^{n+1}\) equals \(1\) for all \(n \).
There is only one power series satisfying this condition, as Hirzebruch \cite{Hirzebruch} observes.

\smallskip
\noindent This series is the \textit{Todd function}\index{generating function!Todd}
(called after John Arthur Todd)%
\begin{equation} {
        \operatorname{T}(z) \,=\, \frac{z}{1- \mathrm{e}^{-z}} \, =
        \,   1 + \frac{1}{2}\frac{z}{1!} + \frac{1}{6}\frac{z^2}{2!} - \frac{1}{30} \frac{z^4}{4!} + \ldots }.
\end{equation}
It generates the Bernoulli numbers matching the values of
the Bernoulli function at the nonnegative integers.

\smallskip
A modern exposition based on the Todd series is the monograph of
Arakawa, Ibukiyama, and Kaneko \cite{Arakawa}. The authors adopt this
definition ``because it is the original definition of Seki and Bernoulli
for one thing, and it is better suited to the special values of the
Riemann zeta function for another.''

\smallskip
Similarly, Neukirch in \textit{Algebraic Number Theory} \cite{Neukirch} calls
the definition \( f(z) = z/(1- \mathrm{e}^{-z}) \) ``more natural and better
suited for the further development of the theory.''
One might hope that all mathematicians will support this simple step towards
greater consistency someday.

\sect{Acknowledgments}{20}

\noindent The author thanks Jörg Arndt, Petros Hadjicostas,
Václav Kotěšovec, Richard J. Mathar, Gergő Nemes, and Michael Somos
for their reading and providing valuable feedback
on an earlier version of this manuscript.
Thanks to Michel Marcus and Jon E. Schoenfield for help with proofreading.

Without using Neil Sloane's \OEIS this essay could have been
written, but it would only have been half as much fun.

\bigskip\bigskip\bigskip
\begin{figure}[h]
    \centering
    \includegraphics[height=90pt,width=130pt]{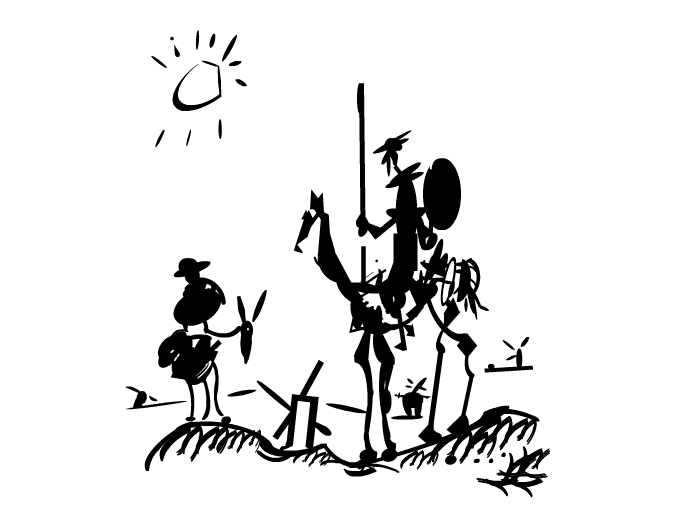}
     \caption[Pablo Picasso, Don Quixote]{\textit{`Don Quixote', sketch by Pablo Picasso.}}
\end{figure}
\clearpage \newpage \FloatBarrier

Some applications of the Swiss-knife polynomials, where

 \( \alpha_n =  {n}/(4^n-2^n), \)  and by convention  \(  \operatorname{K}_{-1}(x) = 1.  \)

\bigskip
{  \centering{  \setlength{\extrarowheight}{3pt}
    \begin{tabular}{cll}
        $\frac{\seqnum{A099612}}{\seqnum{A099617}}$ & \text{Euler zeta}        & { $  \mathcal{E}(n) \, = \,  \lvert \operatorname{K}_n( n \bmod 2) \rvert \,/\,n! $} \label{th3}       \\
        \seqnum{A000111}                            & \text{Euler numbers}     & { $  \operatorname{E}(n) \, = \,   \operatorname{K}_n( n \bmod 2)  $} \label{th6}                      \\
        \seqnum{A028296}                            & \text{Euler secant}      & { $  \operatorname{E}_{\sigma}(2n) \, = \,  \operatorname{K}_{2n}(0) $}  \label{th4}                   \\
        \seqnum{A000182}                            & \text{Euler tangent}     & { $  \operatorname{E}_{\tau}(2n+1) \, = \,  \operatorname{K}_{2n+1}(1) $} \label{th5}                  \\
        $\frac{\seqnum{A164555}}{\seqnum{A027642}}$ & \text{Bernoulli numbers} & { $  \operatorname{B}(n) \, = \, \operatorname{K}_{n-1}(1) \, \alpha_n $}  \label{th7}                 \\
        $\frac{\seqnum{A000367}}{\seqnum{A002445}}$ & \text{Bernoulli tangent} & { $  \operatorname{B}_{\tau}(2n) \, = \, \operatorname{K}_{2n-1}(1) \, \alpha_{2n}  $}  \label{th8}    \\
        $\frac{\seqnum{A160143}}{\seqnum{A193476}}$ & \text{Bernoulli secant}  & { $  \operatorname{B}_{\sigma}(2n+1) \, = \, \operatorname{K}_{2n}(0) \, \alpha_{2n+1} $}  \label{th9} \\
        $\frac{\seqnum{A193472}}{\seqnum{A193473}}$ & \text{Bernoulli extend.} & { $  \mathcal{B}(n) \, = \, \operatorname{K}_{n-1}((n-1)\, \bmod\, 2)\, \alpha_n $} \label{th10}       \\
        \seqnum{A226158}                            & \text{Genocchi}          & { $ \operatorname{G}(n) \, = \, - \operatorname{K}_{n-1}( -1)\,  n \,/\, 2^{n-1} $} \label{th11}       \\
        \seqnum{A188458}                            & \text{Springer}          & { $ \operatorname{S}(n) \, = \, \operatorname{K}_n(1/2)\, 2^n  $}  \label{th12}
    \end{tabular} }   \captionsetup{hypcap=false}}

\bigskip

Concerned with sequences \seqnum{A001896}, \seqnum{A001897},
\seqnum{A019692}, \seqnum{A027642}, \seqnum{A028246}, \seqnum{A081658}, \seqnum{A099612}, \seqnum{A099617},
\seqnum{A122045}, \seqnum{A153641}, \seqnum{A155585}, \seqnum{A160143}, \seqnum{A163626}, \seqnum{A163747},
\seqnum{A164555}, \seqnum{A173018}, \seqnum{A193472}, \seqnum{A193473}, \seqnum{A193476}, \seqnum{A226158},
\seqnum{A263634}, \seqnum{A278075}, \seqnum{A301813}, \seqnum{A303638}, \seqnum{A318259}, \seqnum{A333303},
\seqnum{A335263}, \seqnum{A335750}, \seqnum{A335751}, \seqnum{A335947}, \seqnum{A335948},
\seqnum{A335949}, \seqnum{A335953}, \seqnum{A336454}, \seqnum{A336517}, \seqnum{A337966},
\seqnum{A337967},  \seqnum{A342317}, \seqnum{A344917}, \seqnum{A344918}, \seqnum{A346463}, \seqnum{A346464},
 \seqnum{A346832}, \seqnum{A346833}, \seqnum{A346834}, \seqnum{A346835}, \seqnum{A346838} and \seqnum{A344913}.

\bigskip\bigskip
\hrule

\bigskip

\noindent 2020 MSC: Primary 11B68, Secondary 11M35.

\medskip
\noindent {\it Keywords}:
Bernoulli function, Bernoulli constants, Bernoulli numbers, Bernoulli cumulants,
Bernoulli functional equation, Bernoulli central function, extended Bernoulli function,
Bernoulli central polynomials, alternating Bernoulli function, Stieltjes constants,
Riemann zeta function, Hurwitz zeta function, Worpitzky numbers,
Worpitzky transform, Hasse representation, Hadamard product,  Genocchi numbers,
Genocchi function, Euler secant numbers, Euler tangent numbers, Euler zeta numbers,
Euler function, Euler tangent function, Euler secant function, Eulerian numbers,
Andr\'e numbers, Andr\'e function, Seki numbers, Seki function, Swiss-knife polynomials.

\bigskip

\hrule

\bigskip 

\noindent \textit{Software}: Code repository and supplements, including a Maple worksheet and
a Mathematica Jupyter notebook, available at GitHub \cite{LuschnyIBF}.

\medskip
\noindent \textit{Author}: {\footnotesize ORCID 0000-0001-6245-708X}

\newpage \FloatBarrier

{\Large{\sect{Bernoulli constants}{0}}}

 \bigskip
\( \qquad\qquad\qquad   \beta_n  \quad\qquad\qquad\qquad\qquad\qquad b_n = \frac{\beta_n}{n!}  \)

{\scriptsize {\linespread{1.16}   \begin{Verbatim}
   [ 0]  +1.00000000000000000000000e+00   +1.00000000000000000000000e+00
   [ 1]  -5.77215664901532860606512e-01   -5.77215664901532860606512e-01
   [ 2]  +1.45631690967353449721173e-01   +7.28158454836767248605864e-02
   [ 3]  +2.90710895786169554535912e-02   +4.84518159643615924226519e-03
   [ 4]  -8.21533768121338346464019e-03   -3.42305736717224311026674e-04
   [ 5]  -1.16268503273365002873409e-02   -9.68904193944708357278404e-05
   [ 6]  -4.75994290380637621052001e-03   -6.61103181084218918127779e-06
   [ 7]  +1.67138541801139726910695e-03   +3.31624090875277235933919e-07
   [ 8]  +4.21831653646200836859278e-03   +1.04620945844791874221051e-07
   [ 9]  +3.16911018422735558641847e-03   +8.73321810027379736116201e-09
   [10]  +3.43947744180880481779146e-04   +9.47827778276235895555407e-11
   [11]  -2.25866096399971274152095e-03   -5.65842192760870796637242e-11
   [12]  -3.24221327452684232007482e-03   -6.76868986351369665586675e-12
   [13]  -2.17454785736682251359552e-03   -3.49211593667203185445522e-13
   [14]  +3.84493292452642224040106e-04   +4.41042474175775338023724e-15
   [15]  +3.13813893088949918755710e-03   +2.39978622177099917550506e-15
   [16]  +4.53549848512386314628695e-03   +2.16773122007268285496389e-16
   [17]  +3.39484659125248617003234e-03   +9.54446607636696517342499e-18
   [18]  -4.72986667978530060590399e-04   -7.38767666053863649781558e-20
   [19]  -5.83999975483580370526234e-03   -4.80085078248806522761766e-20
   [20]  -1.00721090609471125811119e-02   -4.13995673771330564126948e-21
   [21]  -9.79321479174274843741249e-03   -1.91682015939912339496482e-22
   [22]  -2.29763093463200254783750e-03   -2.04415431222621660772759e-24
   [23]  +1.24567903906919471380695e-02   +4.81849850110735344392922e-25
   [24]  +2.98550901697978987031938e-02   +4.81185705151256647946111e-26
   [25]  +3.97127819725890390476549e-02   +2.56026331031881493660913e-27
   [26]  +2.79393907712007094428316e-02   +6.92784089530466712388013e-29
   [27]  -1.77336950032031696506289e-02   -1.62860755048558674407104e-30
   [28]  -9.73794335813190698522061e-02   -3.19393756115325557604211e-31
   [29]  -1.85601987419318254285110e-01   -2.09915158936342552768549e-32
   [30]  -2.21134553114167174032372e-01   -8.33674529544144047562508e-34
   [31]  -1.10289594522767989385320e-01   -1.34125937721921866750473e-35
   [32]  +2.40426431930087325860325e-01   +9.13714389129817199794565e-37
   [33]  +8.48223060577873259185100e-01   +9.76842144689316562821221e-38
   [34]  +1.53362895967472747676942e+00   +5.19464288745573322360277e-39
   [35]  +1.78944248075279625487765e+00   +1.73174959516100441594633e-40
   [36]  +7.33429569739007257407068e-01   +1.97162023326628724184554e-42
   [37]  -2.68183987622201934815062e+00   -1.94848008275558832944550e-43
   [38]  -8.96900252642345710339731e+00   -1.71484028164349789185891e-44
   [39]  -1.67295744090075569567376e+01   -8.20162325795024844398325e-46
   [40]  -2.07168737077169487591557e+01   -2.53909617003982347034561e-47
   [41]  -1.01975839351792340684642e+01   -3.04837480681247325379474e-49
   [42]  +3.02221435698546147289327e+01   +2.15103288078139524274228e-50
   [43]  +1.13468181875436585038896e+02   +1.87813767783170614584043e-51
   [44]  +2.31656933743621048467685e+02   +8.71456987091534575015368e-53
   [45]  +3.23493565027658727705394e+02   +2.70429325494165277375631e-54
   [46]  +2.33327851136153134645751e+02   +4.24030297482975157602467e-56
   [47]  -3.10666033627557393250479e+02   -1.20123014394335413203436e-57
   [48]  -1.63390919914361991588774e+03   -1.31619164554817482023699e-58
   [49]  -3.85544150839888666284589e+03   -6.33824832920169119558284e-60
   [50]  -6.29221938159892345466820e+03   -2.06884990450560309799972e-61
 \end{Verbatim}
}}

\clearpage \newpage \FloatBarrier
\begin{table}  {
     \sect{Table of \(\mathcal{S}_n\) and \( \mathcal{Z}_n. \)}{0}}
    \caption{\hspace{-12pt}\textbf{\textcolor{MidnightBlue}{Seki numbers and Euler zeta numbers}}}
 \setlength{\extrarowheight}{4pt}
  \begin{tabular}{c c c}    
        \emph{n} & \(\mathcal{S}_n\) & \( \mathcal{Z}_n \)   \\ 
        \rowcolor{lightgray!20} 0 & $ 1 $ & $ 1 $\\
        \rowcolor{white!50}     1 & $ {\frac{1}{2}} $ & $ 1 $\\
        \rowcolor{lightgray!20} 2 & $ {\frac{1}{6}} $ & $ {\frac{1}{2}} $\\
        \rowcolor{white!50}     3 & $ {\frac{3}{56}} $ & ${\frac{1}{3}} $\\
        \rowcolor{lightgray!20} 4 & $ {\frac{1}{30}} $ & $ {\frac{5}{24}} $\\
        \rowcolor{white!50}     5 & $ {\frac{25}{992}} $ & $ {\frac{2}{15}} $\\
        \rowcolor{lightgray!20} 6 & $ {\frac{1}{42}} $ & $ {\frac{61}{720}} $\\
        \rowcolor{white!50}     7 & $ {\frac{427}{16256}} $ & $ {\frac{17}{315}} $\\
        \rowcolor{lightgray!20} 8 & $ {\frac{1}{30}} $ & $ {\frac{277}{8064}} $\\
        \rowcolor{white!50}     9 & $ {\frac{12465}{261632}} $ & $ {\frac{62}{2835}} $\\
        \rowcolor{lightgray!20} 10 & $ {\frac{5}{66}} $ & $ {\frac{50521}{3628800}} $\\
        \rowcolor{white!50}     11 & $ {\frac{555731}{4192256}} $ & $ {\frac{1382}{155925}} $\\
        \rowcolor{lightgray!20} 12 & $ {\frac{691}{2730}} $ & $ {\frac{540553}{95800320}} $\\
        \rowcolor{white!50}     13 & $ {\frac{35135945}{67100672}} $ & ${\frac{21844}{6081075}} $\\
        \rowcolor{lightgray!20} 14 & $ {\frac{7}{6}} $ & $ {\frac{199360981}{87178291200}} $\\
        \rowcolor{white!50}     15 & $ {\frac{2990414715}{1073709056}} $ & $ {\frac{929569}{638512875}} $\\
        \rowcolor{lightgray!20} 16 & $ {\frac{3617}{510}} $ & $ {\frac{3878302429}{4184557977600}} $\\
        \rowcolor{white!50}     17 & $ {\frac{329655706465}{17179738112}} $ & $ {\frac{6404582}{10854718875}} $\\
        \rowcolor{lightgray!20} 18 & $ {\frac{43867}{798}} $ & $ {\frac{2404879675441}{6402373705728000}} $\\
        \rowcolor{white!50}     19 & $ {\frac{45692713833379}{274877382656}} $ & $ {\frac{443861162}{1856156927625}} $\\
        \rowcolor{lightgray!20} 20 & $ {\frac{174611}{330}} $ & $ {\frac{14814847529501}{97316080327065600}} $\\
        \rowcolor{white!50}     21 & $ {\frac{1111113564712575}{628292059136}} $ & $ {\frac{18888466084}{194896477400625}} $\\
        \rowcolor{lightgray!20} 22 & $ {\frac{854513}{138}} $ & $ {\frac{69348874393137901}{1124000727777607680000}} $\\
        \rowcolor{white!50}     23 & $ {\frac{1595024111042171723}{70368735789056}} $ & $ {\frac{113927491862}{2900518163668125}} $\\
        \rowcolor{lightgray!20} 24 & $ {\frac{236364091}{2730}} $ & $ {\frac{238685140977801337}{9545360026665222144000}} $\\
        \rowcolor{white!50}     25 & $ {\frac{387863354088927172625}{1125899873288192}} $ & $ {\frac{58870668456604}{3698160658676859375}} $\\
        \rowcolor{lightgray!20} 26 & $ {\frac{8553103}{6}} $ & $ {\frac{4087072509293123892361}{403291461126605635584000000}} $\\
        \rowcolor{white!50}     27 & $ {\frac{110350957750914345093747}{18014398375264256}} $ & $ {\frac{8374643517010684}{1298054391195577640625}} $\\
        \rowcolor{lightgray!20} 28 & $ {\frac{23749461029}{870}} $ & $ {\frac{13181680435827682794403}{3209350995912777478963200000}} $\\
        \rowcolor{white!50}     29 & $ {\frac{36315529600705266098580265}{288230375614840832}} $ & $ {\frac{689005380505609448}{263505041412702261046875}} $\\
        \rowcolor{lightgray!20} 30 & $ {\frac{8615841276005}{14322}} $ & $ {\frac{441543893249023104553682821}{265252859812191058636308480000000}} $\\
\end{tabular}\end{table}
\clearpage \newpage \FloatBarrier

\begin{figure}
    \centering
    \includegraphics[height=220pt,width=320pt]{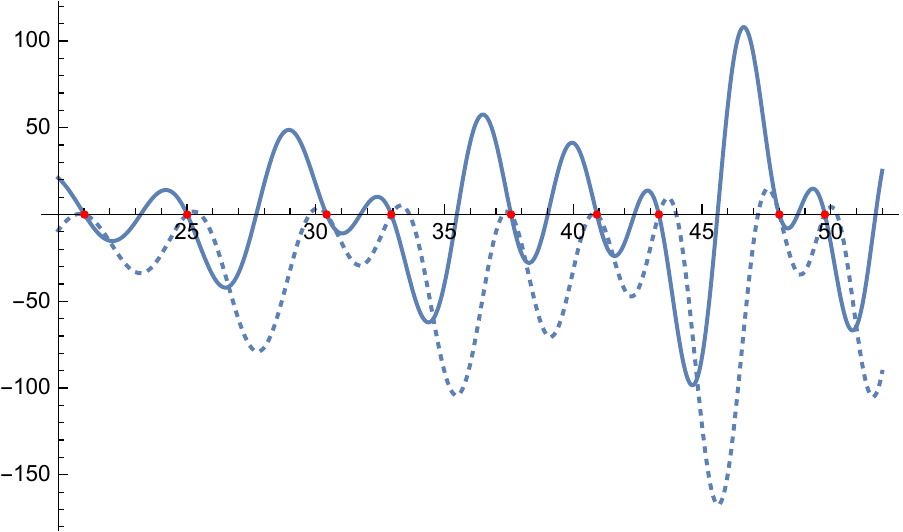}%
    \caption[The Riemann zeros]{\textit{Bernoulli function and Riemann zeros on the critical line.}}
    \label{fig:fig10bernoullizetazeros}
    \label{page-plots}
    \stepcounter{mysection}\addcontentsline{toc}{section}{\arabic{mysection} -- Plots}
\end{figure}

\begin{figure}
    \centering
    \includegraphics[height=240pt,width=320pt]{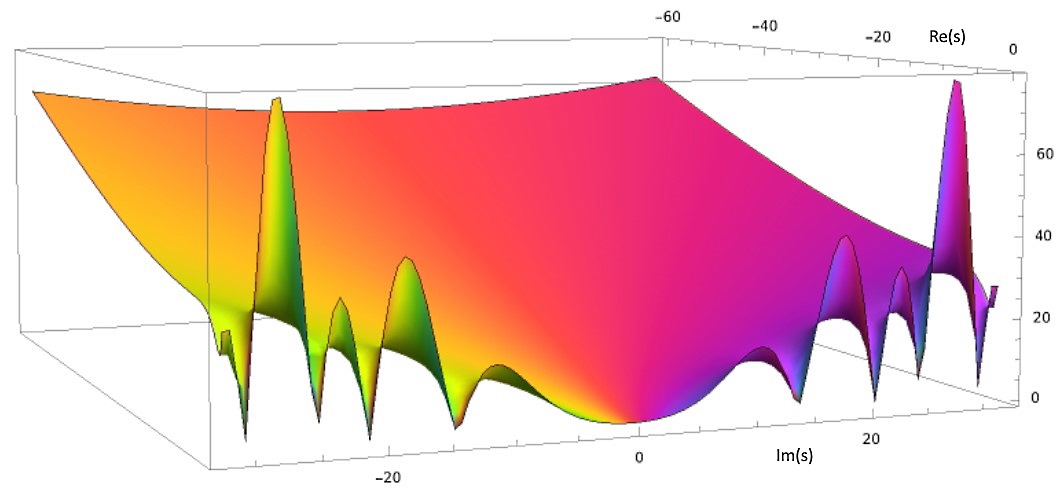}%
    \medskip
    \caption[The Bernoulli tsunami]{\textit{The Bernoulli tsunami.}}
    \vspace{-8pt}
    \caption[]{\textit{The Bernoulli function hits Riemann's critical line.}}
    \label{fig:fig27bernoullitsunami}
\end{figure}

\clearpage \newpage \FloatBarrier

\begin{figure}
    \centering
    \includegraphics[width=320pt]{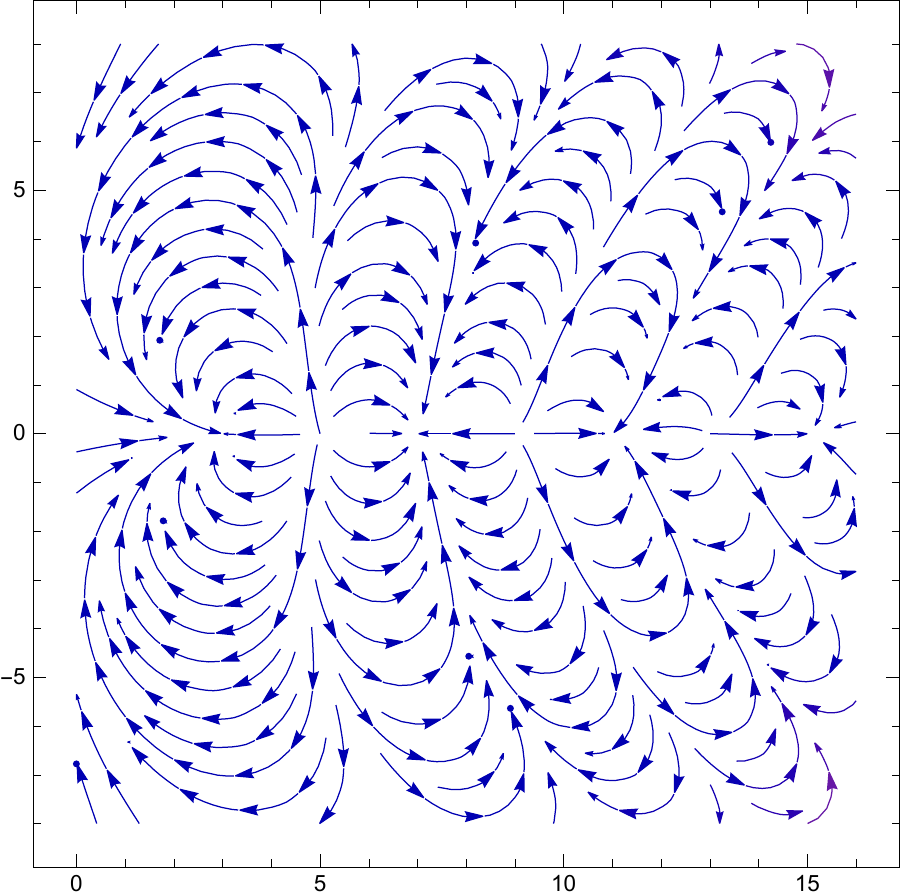}%
    \caption[Phase portrait of \(\operatorname{B}(x)\)]{\ \ \ \ \textit{Phase portrait of the Bernoulli function on the right half plane.}}
    \label{fig:fig11bernoulliphase}
\end{figure}
\clearpage \newpage \FloatBarrier

\begin{figure}
    \centering
    \includegraphics[width=280pt]{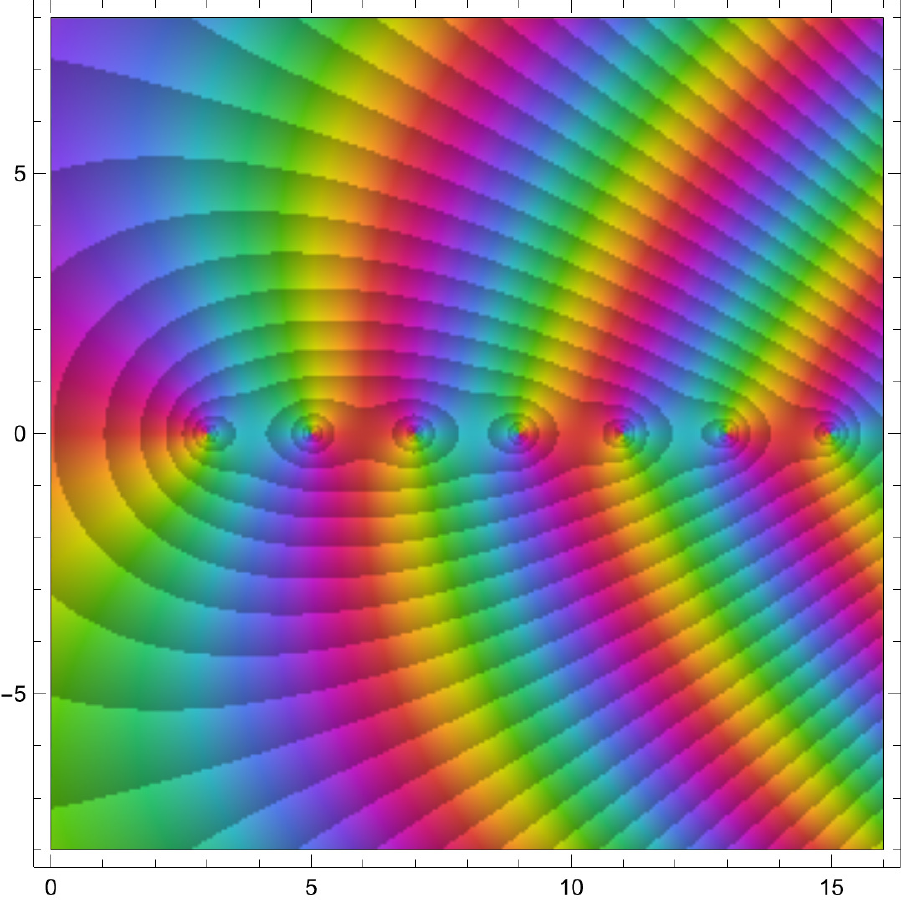}%
    \caption[Complex view of \(\operatorname{B}(x)\)]{\textit{The Bernoulli function on the right half plane, complex view.}}
     \label{fig:fig12bernoullicomplex}
\end{figure}

\begin{figure}
    \centering
    \includegraphics[width=300pt]{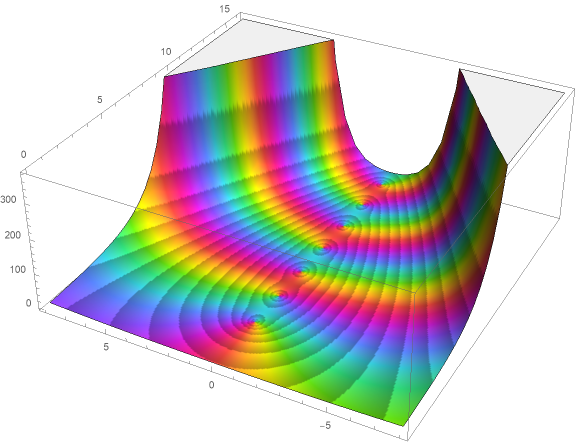}%
    \caption[3-dim view of \(\operatorname{B}(x)\)]{\textit{The Bernoulli function on the right half plane, 3-dim view.}}
    \label{fig:fig13bernoulli3dim}
\end{figure}

\begin{figure}
    \centering
    \includegraphics[width=340pt]{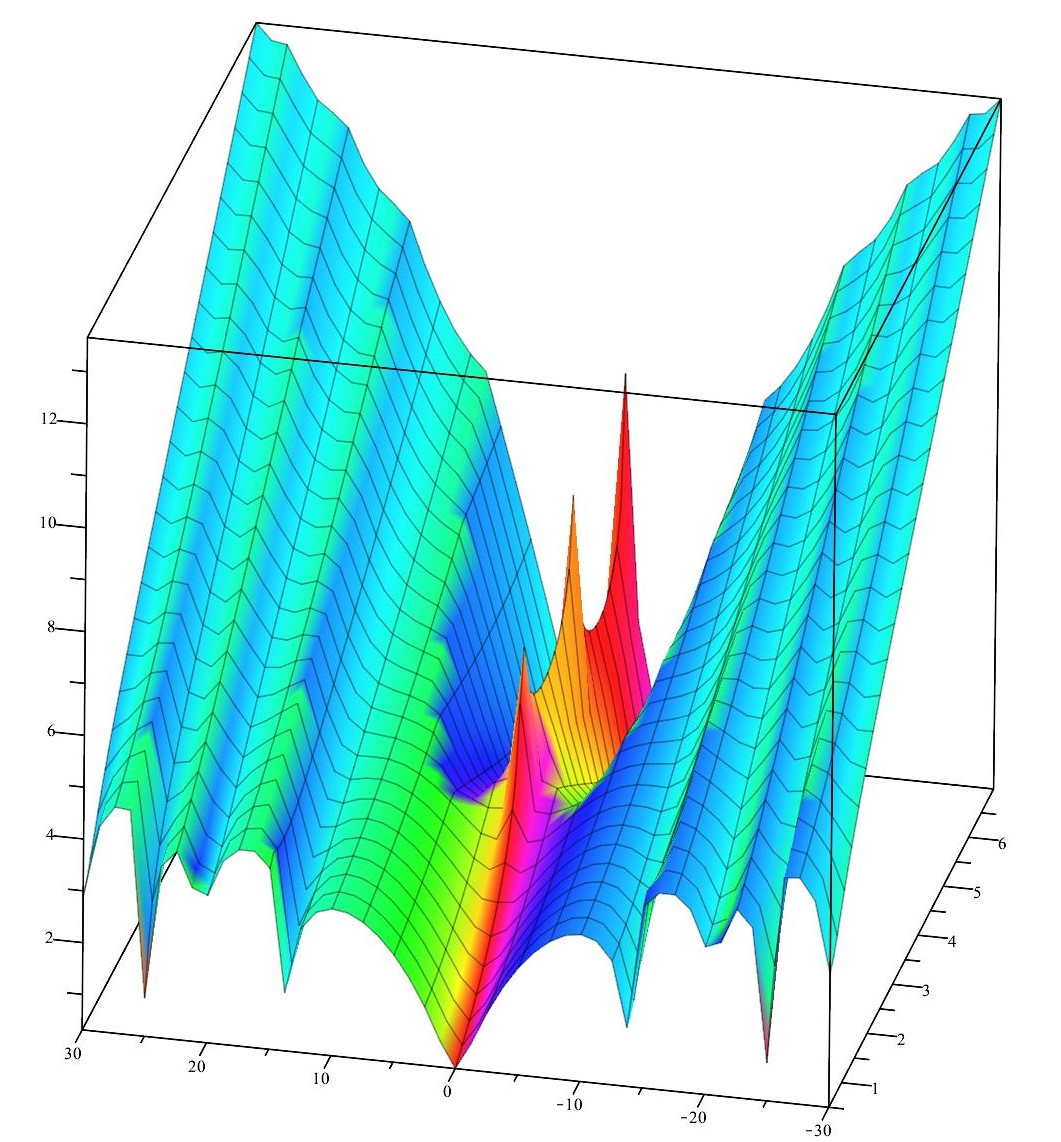}%
    \medskip
    \caption[Logarithm of \(\operatorname{B}(x)\)]{\textit{The logarithm of the Bernoulli function on the right half plane.}}
    \bigskip
\caption[]{The red peaks on the x-axis correspond to the real zero of the Bernoulli
    function (the vanishing of the odd Bernoulli numbers). }
\caption[]{The front side of the plot
    shows the logarithm of the Bernoulli function on the critical line.}
\caption[]{The Hadamard decomposition of the logarithm of the Bernoulli function as seen above is displayed in the two plots below.}
   \label{fig:fig14bernoullilog}
\end{figure}

\begin{figure}
    \centering
    \includegraphics[height=275pt]{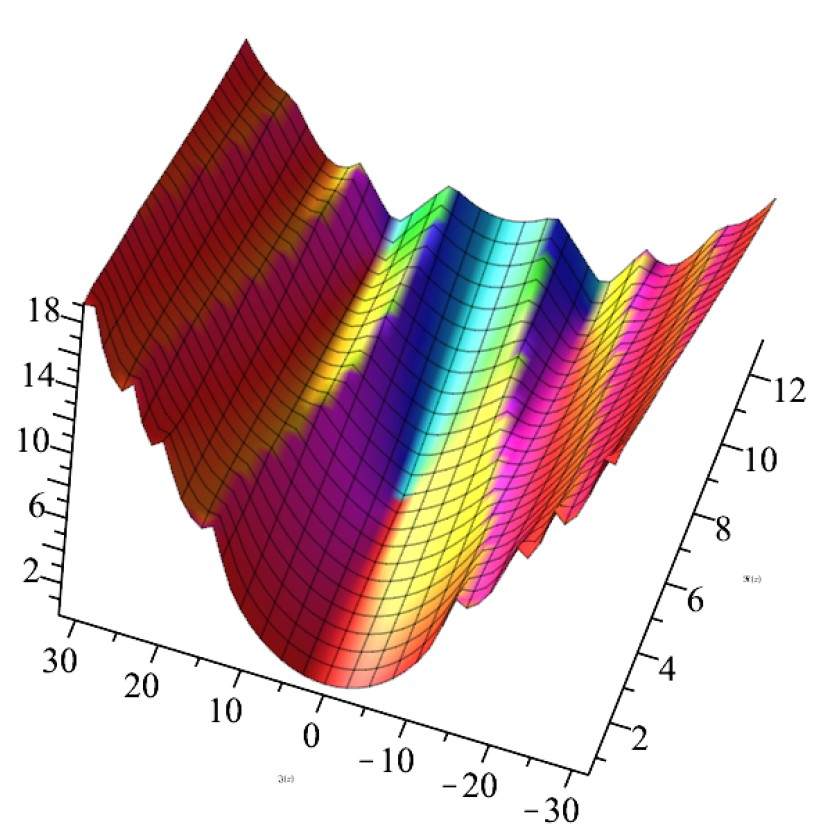}%
    \caption[The Riemann \(\xi\)-factor of \(\operatorname{B}(x)\)]{\textit{The Hadamard decomposition of \(\log \operatorname{B}:\) the Riemann \(\xi\)-factor.}}
    \label{fig:fig15logxi}
    \bigskip
    \includegraphics[height=275pt]{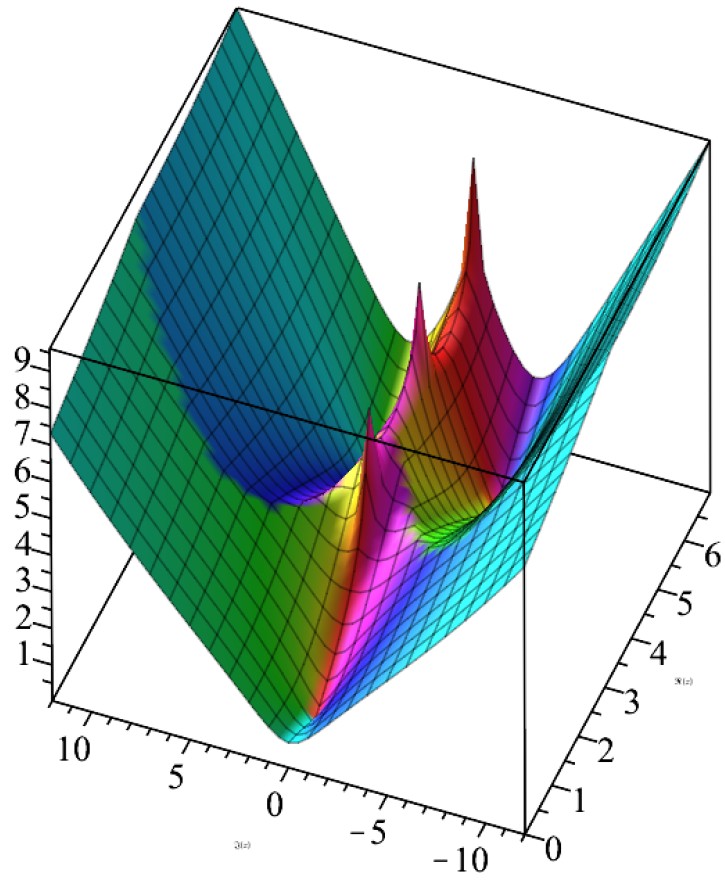}%
    \caption[The singularity factor of \(\operatorname{B}(x)\)]{\textit{The Hadamard decomposition of \(\log \operatorname{B}:\) the singularity factor.}}
    \label{fig:fig16bernoullitrivial}
\end{figure}

\clearpage \newpage \FloatBarrier
\phantomsection\addcontentsline{toc}{section}{References}
\bibliographystyle{crelle}

\newpage
\phantomsection
\label{page-figures}
\addcontentsline{toc}{section}{List of figures}
\listoffigures

\newpage
\phantomsection
\label{page-index}
\addcontentsline{toc}{section}{Index}
\printindex

\cleardoublepage  \newpage \FloatBarrier
\newgeometry{
    left=5cm,
    right=4.2cm,
    top=1.5cm,
    bottom=1.5cm,
    bindingoffset=5mm
}
\thispagestyle{empty}
\phantomsection
\textcolor{MidnightBlue}{
    \vspace{-12pt} }    
    \begingroup
\hypersetup{linkcolor=blue}
\tableofcontents
\endgroup
\label{page-contents}
\end{document}